\documentclass[11pt,a4paper]{article}
\title{\bf Quasi-invariant flow generated by Stratonovich SDE
with BV drift coefficients}
\author{Huaiqian Li$^a$,\quad Dejun Luo$^{b,c}$\footnote{Email: luodj@amss.ac.cn}
\vspace{3mm}\\
{\footnotesize $^a$I.M.B, BP 47870, Universit\'e de Bourgogne,
Dijon, France}\\
{\footnotesize $^b$UR Math\'{e}matiques, Universit\'{e} de
Luxembourg, 6, rue Richard Coudenhove-Kalergi, L-1359 Luxembourg}\\
{\footnotesize $^c$Key Lab of Random Complex Structures and Data
Science, Academy of Mathematics and}\\
{\footnotesize  Systems Science, Chinese Academy of Sciences,
Beijing 100190, China} }
\date{}
\usepackage{amssymb,amsmath,amsfonts,amsthm,color,mathrsfs}

\def\R{\mathbb{R}}
\def\E{\mathbb{E}}
\def\P{\mathbb{P}}

\def\Z{\mathbb{Z}}
\def\L{\mathcal{L}}

\def\d{\textup{d}}
\def\ch{{\bf 1}}
\def\BV{\textup{BV}}
\def\det{\textup{det}}
\def\tr{\textup{tr}}
\def\supp{\textup{supp}}
\def\<{\langle}
\def\>{\rangle}
\def\Proof.{\noindent{\bf Proof. }}
\newcommand{\ra}{\rightarrow}

\newcommand{\ee}{\varepsilon}
\def\F{\mathcal{F}}

\def\div{\textup{div}}
\def\fin{\hfill$\square$}

\def\newdot{{\kern.8pt\cdot\kern.8pt}}
\def\bint{\hskip2pt-\hskip-12pt\int}

\newtheorem{theorem}{Theorem}[section]
\newtheorem{lemma}[theorem]{Lemma}
\newtheorem{corollary}[theorem]{Corollary}
\newtheorem{proposition}[theorem]{Proposition}

\newtheorem{definition}[theorem]{Definition}
\theoremstyle{definition}\newtheorem{remark}[theorem]{Remark}

\begin{document}

\maketitle
\makeatletter 
\renewcommand\theequation{\thesection.\arabic{equation}}
\@addtoreset{equation}{section}
\makeatother 

\begin{abstract}
We generalize the results of Ambrosio [Invent. Math. 158 (2004),
227--260] on the existence, uniqueness and stability of regular
Lagrangian flows of ordinary differential equations to Stratonovich
stochastic differential equations with BV drift coefficients. Then
we construct an explicit solution to the corresponding stochastic
transport equation in terms of the stochastic flow. The approximate
differentiability of the flow is also studied when the drift is a
Sobolev vector field.
\end{abstract}



\section{Introduction}

Let $A_0,A_1,\cdots,A_m$ be vector fields on $\R^d$ and $w_t=(w^1_t,
\cdots, w^m_t)$ an $m$-dimensional standard Brownian motion defined
on a probability space $(\Omega,\F,\P)$. Consider the Stratonovich
stochastic differential equation (abbreviated as SDE)
  \begin{equation}\label{SDE}
  \d X_t=\sum_{i=1}^mA_i(X_t)\circ\d w^i_t+A_0(X_t)\,\d t,\quad X_0=x.
  \end{equation}
It is well known that if $A_i\in C^{2+\delta}_b(\R^d,\R^d)\,
(i=1,\cdots,m)$ and $A_0\in C^{1+\delta}_b(\R^d,\R^d)$, then the
above equation has a unique solution which defines a stochastic flow
of $C^1$-diffeomorphisms on $\R^d$. Here $C_b^{1+\delta}$ means that
the vector field $A_0$ and its first order derivatives are bounded,
and the derivatives are H\"{o}lder continuous of order $\delta>0$.
These conditions on the boundedness of the vector fields and their
derivatives were relaxed in \cite{LiXM}, by allowing the local
Lipschitz constants on the balls centered at the origin to grow as
fast as the logarithmic function. In the case $\delta=0$, it is
proved in \cite{FangImkellerZhang} that under the same growth
conditions, \eqref{SDE} still gives rise to a flow of homeomorphisms
on $\R^d$. This result is generalized in \cite{FangLuo07} to the
case where the drift coefficient $A_0$ satisfies only the general
Osgood condition, at the price of the $C^{3+\delta}_b$ regularity
for the diffusion coefficients; if in addition the distributional
divergence of $A_0$ exists and is bounded, then the Lebesgue measure
is quasi-invariant under the action of the stochastic flow of
homeomorphisms (cf. \cite{Luo09}).

On the other hand, the ordinary differential equation (abbreviated
as ODE)
  \begin{equation}\label{ODE}
  \d X_t=A_0(X_t)\,\d t,\quad X_0=x
  \end{equation}
with Sobolev or even BV coefficient has been studied intensively in
the last three decades. Here $A_0$ can be a time-dependent vector
field. The existence of quasi-invariant flow of measurable maps
associated to a vector field $A_0$ with Sobolev regularity was first
studied by Cruzeiro \cite{Cruzeiro83}. A thorough treatment was
carried out by DiPerna and Lions in the celebrated paper
\cite{DiPernaLions}, where the authors deduced the existence and
uniqueness of flows generated by \eqref{ODE} from the well posedness
of the corresponding transport equations. Similar results were
obtained in \cite{CiprianoCruzeiro05} by taking the standard
Gaussian measure as the reference measure. Ambrosio
\cite{Ambrosio04} generalized the results to the case where $A_0$
has only BV regularity by considering the continuity equation. S.
Fang \cite{Fang09} gave a short introduction to the theories
mentioned above. The extension of these results to the infinite
dimensional Wiener space have been done in \cite{AmbrosioFigalli09,
FangLuo10}. Using the local maximal function, Crippa and De Lellis
obtained in \cite{CrippaLellis} some new estimates which allow them
to give a direct proof of the existence and uniqueness of the
DiPerna-Lions flow.

Inspired by these studies of ODE, there have been several attempts
to solve the SDE with Sobolev coefficients. Following the method in
\cite{CrippaLellis}, X. Zhang \cite{Zhang09} showed the existence
and uniqueness of the stochastic flow of measurable maps generated
by It\^{o} SDE with Sobolev coefficients, provided that the
derivatives of the diffusion vector fields are bounded. The SDE with
BV drift vector field was also considered in this paper, but the
diffusion coefficients were assumed to be constant. In
\cite{FangLuoThalmaier}, the authors took the standard Gaussian
measure as the reference measure and proved a priori estimate on the
$L^p$ norm of the density of the flow, which enabled them to
construct the unique flow associated to \eqref{SDE}, provided that
the gradients of the diffusion coefficients and the divergences with
respect to the Gaussian measure are exponentially integrable. In the
recent work \cite{Zhang10}, X. Zhang studied the Stratonovich SDE
with drift coefficient belonging to $W^{1,1}_{loc}(\R^d,\R^d)$, and
he also provided a Freidlin-Wentzell type large deviation estimate
for the stochastic flow.

In the present paper we consider the Stratonovich SDE \eqref{SDE}
with BV drift vector field. Our method is based on Ocone-Pardoux's
decomposition \cite{OconePardoux} of the flow generated by
\eqref{SDE} into the stochastic flow of the diffusion part, and a
flow associated to random ODEs whose driving vector field is a
transformation of the drift coefficient $A_0$ by the stochastic
flow. This approach was applied in \cite{FangLuo07} to deal with the
Stratonovich SDE with drift satisfying the general Osgood condition.
Our main result can be stated as follows ($\L_d$ is the Lebesgue
measure on $\R^d$).

\begin{theorem}\label{sect-1-thm-1}
Assume that $A_1,\cdots,A_m\in C^{3+\delta}_b(\R^d, \R^d)$, and the
drift $A_0$ satisfies
  \begin{enumerate}
  \item[\rm(1)] $A_0$ has sublinear growth, i.e. $|A_0(x)|\leq
  C(1+|x|^{1-\ee_0}),\,x\in\R^d$;
  \item[\rm(2)] $A_0\in \BV_{loc}(\R^d,\R^d)$;
  \item[\rm(3)] the divergence $D\cdot A_0\ll\L_d$ with locally
  bounded density function $\div(A_0)$.
  \end{enumerate}
Then the equation \eqref{SDE} generates a unique stochastic flow of
measurable maps $X_t:\R^d\ra\R^d$ which leaves the Lebesgue measure
$\L_d$ quasi-invariant.
\end{theorem}

This result will be proved in Section 3. Note that the sublinear
growth of $A_0$ ensures that the vector field $\tilde A_0$ defined
in \eqref{sect-3.1} has similar growth (Lemma \ref{sect-3-lem-2}),
which in turn implies the classical growth estimates on the solution
of the ODE. The sublinear growth of $A_0$ also allows us to assume
only the local boundedness of the divergence $\div(A_0)$ (as in (ii)
of Theorem 6.2 in \cite{Ambrosio04}), compared to the global
boundedness required in \cite{Ambrosio08, DiPernaLions}. The
starting point of the proof of Theorem \ref{sect-1-thm-1} is the
relationship between the distributional derivative (a Radon measure
in the present case) of the drift $A_0$ and that of the transformed
vector field $\tilde A_0$ defined in \eqref{sect-3.1}. This will be
done in Lemma \ref{sect-2-lem-1}, where we also show that if the
divergence $D\cdot A_0$ is absolutely continuous with respect to
$\L_d$, then so is $D\cdot\tilde A_0$, and they are related to each
other by a quite simple equality.

We study in Section 4 the stability of the solutions to \eqref{SDE}
when a sequence of vector fields $A^n_0$ converge in some sense to
$A_0$. For this purpose, denote by $B(R)$ the ball centered at the
origin with radius $R$.

\begin{theorem}\label{sect-1-thm-2}
Assume the conditions of Theorem \ref{sect-1-thm-1}. Let $A^n_0:\R^d
\ra\R^d$ be vector fields satisfying
  \begin{enumerate}
  \item[\rm(1)] there is $C>0$ and $\ee_0\in(0,1)$, such that
  $\sup_{n\geq1}|A^n_0(x)|\leq C(1+|x|^{1-\ee_0}),\,x\in\R^d$;
  \item[\rm(2)] $A^n_0$ converges to $A_0$ in $L^1_{loc}(\R^d,
  \R^d)$;
  \item[\rm(3)] for any $n\geq 1$, $\nabla A^n_0$ is locally
  bounded;
  \item[\rm(4)] for any $R>0$, $\sup_{n\geq1}\|\div(A^n_0)\|_{L^\infty(B(R))}
  <+\infty$.
  \end{enumerate}
Let $X^n_t$ be the flow associated to \eqref{SDE} with $A_0$ being
replaced by $A^n_0$. Then for any $p>1$ and $T,\,R>0$, the following
convergence holds almost surely and in $L^p(\Omega,\P)$:
  $$\lim_{n\ra\infty}\int_{B(R)}\sup_{t\in[0,T]}|X^n_t(x)-X_t(x)|\,\d x=0.$$
\end{theorem}

It is well known that when the coefficients are smooth, the solution
to the corresponding stochastic transport equation can be explicitly
expressed in terms of the flow generated by \eqref{SDE}, see
\cite{FangLuo07} Theorem 5.1 for the case where $A_0$ satisfies the
general Osgood condition. As an application of the above stability
result, we will show that similar representation still holds even in
the situation of BV drift.

Finally in Section 5 we consider slightly more regular drift
coefficient $A_0\in W^{1,1}_{loc}$, showing that almost surely, the
stochastic flow $X_t$ is approximately differentiable on $\R^d$.
This generalizes the results in \cite{AmbrosioLecumberryManiglia,
CrippaLellis} to the stochastic context.

\section{Preparations and known results of ODE}

In this section we give some preliminary results needed in the
subsequent sections. Here is the definition of BV functions.

\begin{definition}\label{sect-2-def-1}
A locally integrable function $b:\R^d\ra\R^m$ is said to be of class
$\BV_{loc}$ if there is a $\R^{m\times d}$-valued Radon measure
$Db=(D_ib^j)_{1\leq i\leq d,\,1\leq j\leq m}$ such that
  $$\int_{\R^d}\psi\,\d(D_ib^j)=-\int_{\R^d}b^j\partial_i\psi\,\d x,
  \quad \mbox{for all }\psi\in C_c^1(\R^d).$$
\end{definition}

If $b\in C^1$, then we still denote by $Db$ the function $\nabla b$
on $\R^d$. If $m=d$, we denote by $D\cdot
b=\tr(Db)=\sum_{i=1}^dD_ib^i$ the divergence of the BV vector field
$b$, which is again a Radon measure on $\R^d$. In the following
$\det(\cdot)$ means the determinant of a matrix. For a measurable
map $\varphi:\R^d\ra\R^d$ and a Radon measure $\mu$ on $\R^d$,
$\varphi_\#\mu$ denotes the push forward of the measure $\mu$ by
$\varphi$ (or the ``distribution'' of $\varphi$ under $\mu$). Now we
prove

\begin{lemma}\label{sect-2-lem-1}
Let $b:\R^d\ra\R^d$ be a $\BV_{loc}$ vector field and
$\varphi:\R^d\ra\R^d$ a $C^2$-diffeomorphism, then
\begin{enumerate}
\item[\rm(1)] the composition $b\circ\varphi$ is still a $\BV_{loc}$
vector field and
  $$D(b\circ\varphi)=|\det(J_\varphi)|^{-1}J_\varphi^\ast\,[(\varphi^{-1})_\# Db],$$
where $J_\varphi$ is the Jacobi matrix of $\varphi$ and
$J_\varphi^\ast$ is its transpose,
$J_\varphi^\ast\,[(\varphi^{-1})_\# Db]$ is the matrix product of
$J_\varphi^\ast$ and $(\varphi^{-1})_\# Db$;
\item[\rm(2)] if the divergence $D\cdot b\ll \L_d$ with density
function $\div(b)$, then
  $$D\cdot[J_\varphi^{-1}(b\circ\varphi)]=\<\div(J_\varphi^{-1}),
  b\circ\varphi\>+\div(b)\circ\varphi,$$
where $\div(J_\varphi^{-1})$ is a vector field whose components are
the divergences of the column vectors of $J_\varphi^{-1}$.
\end{enumerate}
\end{lemma}

\Proof. (1) For every $i\in\{1,\cdots,d\}$, we only have to show
that $b^i\circ\varphi$ is a BV function. By Theorems 2 and 3 in
Section 5.2 of \cite{EvansGariepy}, there exists a sequence of
functions $b_n^i\in \BV_{loc}\cap C^\infty(\R^d)$ such that
$b_n^i\ra b^i$ in $L^1_{loc}(\R^d)$ and the vector valued measures
$Db_n^i$ converges weakly to $Db^i$ as $n\ra\infty$. Note that the
composition $b_n^i\circ \varphi$ belongs to $C^2$, hence for any
$\psi\in C^\infty_c(\R^d)$, by the integration by parts formula,
  \begin{equation}\label{sect-2-lem-1.1}
  \int_{\R^d}\psi\,[D(b_n^i\circ\varphi)]\,\d x=-\int_{\R^d}(b_n^i\circ\varphi)
  \nabla\psi\,\d x.
  \end{equation}
We have by the chain rule, $D(b_n^i\circ\varphi)
=J_\varphi^\ast\,[(Db_n^i)\circ \varphi]$. It follows from the
formula of changing variables that
  \begin{eqnarray*}
  \int_{\R^d}\psi\,[D(b_n^i\circ\varphi)]\,\d x
  &=&\int_{\R^d}\psi J_\varphi^\ast\,[(Db_n^i)\circ\varphi]\,\d x\cr
  &=&\int_{\R^d}[(\psi J_\varphi^\ast)\circ(\varphi^{-1})](Db_n^i)|\det(J_{\varphi^{-1}})|\,\d
  x\cr
  &\ra&\int_{\R^d}[(\psi J_\varphi^\ast)\circ(\varphi^{-1})]\cdot|\det(J_{\varphi^{-1}})|\,\d(Db^i)
  \end{eqnarray*}
as $n\ra\infty$, due to the weak convergence of $(Db_n^i)\,\d x$ to
$Db^i$. Now we have
  \begin{eqnarray*}
  \int_{\R^d}[(\psi J_\varphi^\ast)\circ(\varphi^{-1})]\cdot|\det(J_{\varphi^{-1}})|\,\d(Db^i)
  &=&\int_{\R^d}\psi |\det(J_{\varphi^{-1}})\circ\varphi| J_\varphi^\ast\,\d[(\varphi^{-1})_\#
  Db^i]\cr
  &=&\int_{\R^d}\psi |\det(J_{\varphi})|^{-1} J_\varphi^\ast\,\d[(\varphi^{-1})_\#
  Db^i],
  \end{eqnarray*}
where the last equality follows from $J_{\varphi^{-1}}\circ \varphi=
(J_\varphi)^{-1}$. Therefore
  \begin{equation}\label{sect-2-lem-1.2}
  \lim_{n\ra\infty}\int_{\R^d}\psi\,[D(b_n^i\circ\varphi)]\,\d x
  =\int_{\R^d}\psi |\det(J_{\varphi})|^{-1} J_\varphi^\ast\,\d[(\varphi^{-1})_\#
  Db^i].
  \end{equation}

Now we consider the limit of the right hand side of
\eqref{sect-2-lem-1.1}. Again by changing variables,
  \begin{eqnarray*}
  -\int_{\R^d}(b_n^i\circ\varphi)\nabla\psi\,\d x
  &=&-\int_{\R^d}b_n^i\cdot[(\nabla\psi)\circ\varphi^{-1}]\cdot|\det(J_{\varphi^{-1}})|\,\d
  x
  \end{eqnarray*}
which converges to
  $$-\int_{\R^d}b^i\cdot[(\nabla\psi)\circ\varphi^{-1}]\cdot|\det(J_{\varphi^{-1}})|\,\d
  x
  =-\int_{\R^d}(b^i\circ\varphi)\nabla\psi\,\d x.$$
This combines with \eqref{sect-2-lem-1.1} and \eqref{sect-2-lem-1.2}
leads to
  $$\int_{\R^d}\psi |\det(J_{\varphi})|^{-1} J_\varphi^\ast\,\d[(\varphi^{-1})_\# Db^i]
  =-\int_{\R^d}(b^i\circ\varphi)\nabla\psi\,\d x,$$
which means that $b^i\circ\varphi$ is a BV function and
  $$D(b^i\circ\varphi)=|\det(J_{\varphi})|^{-1} J_\varphi^\ast\,[(\varphi^{-1})_\# Db^i].$$

(2) By the chain rule and (1), it is easy to know that
$J_\varphi^{-1}(b\circ\varphi)$ is a $\BV_{loc}$ vector field, and
  \begin{equation}\label{sect-2-lem-1.3}
  D\cdot[J_\varphi^{-1}(b\circ\varphi)]
  =\<\div(J_\varphi^{-1}),b\circ\varphi\>
  +\<J_\varphi^{-1},D(b\circ\varphi)\>,
  \end{equation}
where the second $\<\cdot,\cdot\>$ is the inner product of matrices
regarded as elements in $\R^{d\times d}$. By the expression in (1),
we have
  \begin{equation}\label{sect-2-lem-1.4}
  \<J_\varphi^{-1},D(b\circ\varphi)\>=|\det(J_\varphi)|^{-1}\tr[(\varphi^{-1})_\# Db]
  =|\det(J_\varphi)|^{-1}(\varphi^{-1})_\#\tr(Db).
  \end{equation}
Since $\tr(Db)=D\cdot b=\div(b)\L_d$, for any $\psi\in
C^\infty_c(\R^d)$,
  \begin{eqnarray*}
  \int_{\R^d}\psi\,\d[(\varphi^{-1})_\#\tr(Db)]&=&\int_{\R^d}(\psi\circ\varphi^{-1})\,\d[\tr(Db)]
  =\int_{\R^d}(\psi\circ\varphi^{-1})\cdot\div(b)\,\d x\cr
  &=&\int_{\R^d}\psi\cdot[\div(b)\circ\varphi]\cdot|\det(J_\varphi)|\,\d x.
  \end{eqnarray*}
Therefore $(\varphi^{-1})_\#\tr(Db)=[\div(b)\circ\varphi]\cdot
|\det(J_\varphi)|$. Combining this equality with
\eqref{sect-2-lem-1.3} and \eqref{sect-2-lem-1.4}, we complete the
proof. \fin

\medskip

The next technical result will be used in Section 4.

\begin{lemma}\label{sect-2-lem-2}
Let $\varphi:\R^d\ra\R^d$ be a $C^2$-diffeomorphism and $\tilde\rho
:=\frac{\d(\varphi_\#\L_d)}{\d\L_d}$ the Radon-Nikodym density
function. Then we have
  $$(\tilde\rho^{-1}\nabla\tilde\rho)\circ\varphi=\div\big(J_\varphi^{-1}\big).$$
\end{lemma}

\Proof. It is well known that $\tilde\rho\circ\varphi=
\big|\det\big(J_\varphi^{-1} \big)\big|$. Since $\varphi$ is a
$C^2$-diffeomorphism of $\R^d$, the function $\R^d\ni
x\mapsto\det\big(J_\varphi^{-1} \big)(x)$ does not change sign.
Without loss of generality, we may assume that
$\det\big(J_\varphi^{-1} \big)>0$ on the whole $\R^d$, thus
$\tilde\rho\circ\varphi= \det\big(J_\varphi^{-1} \big)$. As
$\nabla(\tilde\rho\circ \varphi)=J_{\varphi}^\ast
[(\nabla\tilde\rho)\circ\varphi]$, we have
  $$(\nabla\tilde\rho)\circ\varphi=(J_{\varphi}^{-1})^\ast\nabla(\tilde\rho\circ
  \varphi)=(J_{\varphi}^{-1})^\ast\nabla\det\big(J_\varphi^{-1}\big).$$
Note that $\nabla\det\big(J_\varphi^{-1}\big)
=-[\det(J_\varphi)]^{-2}\nabla\det(J_\varphi)$, so the equality that
we should prove is
  $$\nabla\det(J_\varphi)=-\det(J_\varphi)
  J_{\varphi}^\ast\,\div\big(J_\varphi^{-1}\big).$$
For simplification of the notations we write $J=J_\varphi$ and
$K=J_\varphi^{-1}$, then $J_{ij}=\partial_j\varphi^i,\,1\leq i,j\leq
d$. We consider the determinant $\det(\cdot)$ as a smooth function
defined on $\R^{d\times d}$. By the chain rule and Jacobi's formula
(see \cite{Magnus} Part Three, Section 8.3), we have for any
$l\in\{1,\cdots,d\}$,
  \begin{equation}\label{sect-2-lem-2.1}
  \partial_l\det(J)=\sum_{i,j=1}^d\frac{\partial\det}{\partial x_{ij}}(J)
  \cdot\partial_lJ_{ij}=\sum_{i,j=1}^d\det(J)K_{ji}\partial_{lj}\varphi^i.
  \end{equation}
For any $1\leq j\leq d$, it holds $\delta_{jl}=\sum_{i=1}^dK_{ji}
J_{il}=\sum_{i=1}^dK_{ji}\partial_{l}\varphi^i$, therefore
  $$\sum_{i=1}^d(\partial_{j}K_{ji})(\partial_l\varphi^i)
  +\sum_{i=1}^dK_{ji}(\partial_{jl}\varphi^i)=0,$$
that is, $\sum_{i=1}^dK_{ji}(\partial_{jl}\varphi^i)=-\sum_{i=1}^d
(\partial_{j}K_{ji})(\partial_l\varphi^i)$. Combining this with
\eqref{sect-2-lem-2.1} leads to
  \begin{eqnarray*}
  \partial_l\det(J)&=&-\det(J)\sum_{j=1}^d\sum_{i=1}^d(\partial_{j}K_{ji})(\partial_l\varphi^i)
  =-\det(J)\sum_{i=1}^d(\partial_l\varphi^i)\sum_{j=1}^d(\partial_{j}K_{ji})\cr
  &=&-\det(J)\sum_{i=1}^dJ_{il}\div(K_{\cdot
  i})=-\det(J)(J^\ast\div(K))_l.
  \end{eqnarray*}
The proof is complete. \fin

\medskip

Now we recall the definition of the regular Lagrangian flow
associated to a time-dependent vector field $b_t$ (see
\cite{Ambrosio08, CrippaLellis}).

\begin{definition}\label{sect-2-def-2} Let $b\in L^1_{loc}([0,T]\times \R^d,\R^d)$. We
call a map $Y:[0,T]\times \R^d\ra\R^d$ a regular Lagrangian flow for
the vector field $b$ if
\begin{enumerate}
\item[\rm(1)] for a.e. $x\in\R^d$, the map $[0,T]\ni t\ra Y_t(x)$ is an
absolutely continuous integral solution of
  \begin{equation}
  \d Y_t=b_t(Y_t)\,\d t,\quad Y_0=x;
  \end{equation}
\item[\rm(2)] $(Y_t)_\#\L_d \ll\L_d$ for all $t\in[0,T]$.
\end{enumerate}
\end{definition}

Note that this definition is slightly different from that in
\cite{CrippaLellis}: in condition (2) we do not require that
$(Y_t)_\#\L_d $ is dominated by $C\L_d$ on the whole $\R^d$, where
$C>0$ is a constant.

Given a measurable map $Y:\R^d\ra\R^d$, we say that $Z:\R^d\ra\R^d$
is a measurable inverse map of $Y$ if $Z$ is measurable and for a.e.
$x\in\R^d$, $x=Y(Z(x))=Z(Y(x))$. We often denote by $Y^{-1}$ the
measurable inverse map of $Y$ (see \cite{Zhang09} Lemma 3.4 for a
characterization of this notion). In the following theorem we
summarize the results concerning the existence and uniqueness of
regular Lagrangian flow generated by a BV vector field.

\begin{theorem}\label{sect-2-thm-2}
Let $b_t:\R^d\ra\R^d$ be a time-dependent vector field satisfying:
\begin{enumerate}
\item[\rm(1)] $\frac{|b_t(x)|}{1+|x|}\in L^\infty([0,T]\times\R^d)$;
\item[\rm(2)] $b_t\in\BV_{loc}(\R^d)$ for a.e. $t\in[0,T]$, and for
any $R>0$, $|Db_t|(B(R))\in L^1_{loc}(0,T)$ and
  $$\int_0^T\|\div(b_t)\|_{L^\infty(B(R))}\d t<+\infty.$$
\end{enumerate}
Then the vector field $b$ generates a unique regular Lagrangian flow
$\{Y_t:0\leq t\leq T\}$. Moreover for any $t\in[0,T]$, $Y_t$ has a
measurable inverse map $Y_t^{-1}:\R^d \ra\R^d$ and
$(Y_t^{-1})_\#\L_d=\rho_t\,\L_d$ with
  $$\rho_t(x)=\exp\bigg(\int_0^t\div(b_s)(Y_s(x))\,\d s\bigg).$$
\end{theorem}

\Proof. The first part of this theorem was first proved in
\cite{Ambrosio04} for bounded vector field $b$, and then in
\cite{Ambrosio08} for the general case (see the remark at the end of
Section 5), while the second part was proved in \cite{Zhang09} for a
BV vector field $b$ independent of time (just let the diffusion
coefficients $\sigma$ be 0 in Theorem 2.6), but the proof for the
time-dependent case is similar. \fin

\begin{remark}\label{sect-2-rem-1}
For the density function $\tilde \rho_t:=\frac{\d((Y_t)_\#\L_d)}{
\d\L_d}$, we have
  $$\tilde\rho_t(x)=\exp\bigg(-\int_0^t\div(b_s)\big[Y_s(Y_t^{-1}(x))\big]\,\d s\bigg).$$
See \cite{CiprianoCruzeiro05} Theorem 2.1 where the expressions are
given using double time parameters.
\end{remark}

\section{Existence and uniqueness of \eqref{SDE} with BV drift}

In this section we prove Theorem \ref{sect-1-thm-1}. First we
introduce Ocone and Pardoux's decomposition of the Stratonovich SDE
\eqref{SDE} (see \cite{OconePardoux} PART II or \cite{FangLuo07}
Section 2). Consider the following Stratonovich SDE without drift:
  $$\d\tilde X_t=\sum^m_{i=1}A_i(\tilde X_t)\circ\d w^i_t,\quad \tilde X_0=x. $$
It is well-known that under the conditions that $A_1, \cdots, A_m\in
C_b^{3+\delta}$ for some $\delta>0$, the solutions of the above SDE
admit a version $\tilde X_t(x,w)$ such that there exists a full
subset $\Omega_0$, for each $w\in\Omega_0$ and each $t>0$, $x\ra
\tilde X_t(x,w)$ is a $C^2$-diffeomorphism of $\R^d$. Set
$\varphi_t(x)=\tilde X_t(x,w)$. Let $J_t(x)=\partial_x\varphi_t(x)$
be the Jacobian matrix of $\varphi_t:\R^d\ra\R^d$ and
$K_t(x)=(J_t(x))^{-1}$ the inverse of $J_t(x)$. Define for $w\in
\Omega_0$,
  \begin{equation}\label{sect-3.1}
  \tilde A_0(t,x)=K_t(x)A_0(\varphi_t(x)).
  \end{equation}
We consider the differential equation
  \begin{equation}\label{sect-3.2}
  \d Y_t=\tilde A_0(t,Y_t)\,\d t,\quad Y_0=x.
  \end{equation}
Then the solutions of \eqref{SDE} can be expressed as (at least when
$A_0$ is smooth)
  \begin{equation}\label{sect-3.3}
  X_t(x)=\varphi_t(Y_t(x)).
  \end{equation}
Therefore it is sufficient to study the well posedness of the random
ODE \eqref{sect-3.2} under the assumptions on $A_0$ in Theorem
\ref{sect-1-thm-1}, and then show that the representation
\eqref{sect-3.3} indeed gives the flow associated to the original
Stratonovich SDE \eqref{SDE}.

In the next lemma, we collect various growth results concerning the
stochastic flow $\varphi_t$ and its derivatives, and the random
vector field $\tilde A_0$ defined in \eqref{sect-3.1} (see
\cite{FangLuo07} Lemma 2.2 for a proof).

\begin{lemma}\label{sect-3-lem-2}
Assume the conditions of Theorem \ref{sect-1-thm-1}, we have for any
$T>0$,
\begin{enumerate}
\item[\rm(1)] for any $\alpha>1$ and $\beta>0$, there is $F,G\in
\cap_{p>1}L^p(\Omega)$ such that for all $(t,x)\in[0,T]\times \R^d$,
  $$|\varphi_t(x)|\leq F\cdot(1+|x|^\alpha),
  \quad\|J_t(x)\|\vee\|K_t(x)\|\leq G\cdot(1+|x|^\beta);$$
\item[\rm(2)] there exist $\ee_1\in(0,1)$ and $\Phi_T\in
\cap_{p>1}L^p(\Omega)$, such that for all $(t,x)\in[0,T]\times
\R^d$,
  $$|\tilde A_0(t,x)|\leq \Phi_T(1+|x|^{1-\ee_1}).$$
\end{enumerate}
\end{lemma}

Now we can prove

\begin{proposition}\label{sect-3-prop-1}
Under the conditions of Theorem \ref{sect-1-thm-1}, for almost
surely $w\in\Omega_0$, the ODE \eqref{sect-3.2} generates a unique
regular Lagrangian flow $Y_t$ which leaves the Lebesgue measure
quasi-invariant; moreover for all $t\geq0$, $Y_t$ has a measurable
inverse map $Y_t^{-1}:\R^d\ra\R^d$.
\end{proposition}

\Proof. We only need to check that the conditions in Theorem
\ref{sect-2-thm-2} are satisfied for almost surely $w\in \Omega_0$.
First by (2) of Lemma \ref{sect-3-lem-2}, the vector field $\tilde
A_0$ satisfies the condition (1) in Theorem \ref{sect-2-thm-2}.

For all $t\in[0,T]$, since $\varphi_t$ is a $C^2$-diffeomorphism on
$\R^d$ and $K_t:\R^d\ra\R^{d\times d}$ is $C^1$, Lemma
\ref{sect-2-lem-1} and the definition \eqref{sect-3.1} of $\tilde
A_0$ tell us that $\tilde A_0(t,\cdot)\in \BV_{loc}(\R^d,\R^d)$.
Moreover by the chain rule and Lemma \ref{sect-2-lem-1} (1),
  \begin{eqnarray}\label{sect-3-prop-1.1}
  D\tilde A_0(t)&=&(DK_t)A_0(\varphi_t)+K_tD(A_0(\varphi_t))\cr
  &=&(DK_t)A_0(\varphi_t)+|\det(J_t)|^{-1}K_t\,[(\varphi_t^{-1})_\# DA_0]J_t.
  \end{eqnarray}
Fix any $R>0$. Recall that we identify a locally integrable function
$f$ on $\R^d$ with the Radon measure $f\,\d x$. Since $K_t\in C^1$,
we know that $(t,x)\ra DK_t(x)$ is bounded on $[0,T]\times B(R)$.
Moreover, by the sublinear growth of $A_0$ and Lemma
\ref{sect-3-lem-2} (1), it is easy to deduce the boundedness of
$A_0(\varphi_t)$ on $[0,T]\times B(R)$. Hence there exists a
positive constant $C_{T,R}$ (depends on $w\in\Omega_0$) such that
the total variation
  $$|(DK_t)A_0(\varphi_t)|(B(R))\leq C_{T,R} \L_d(B(R)),$$
which implies that $t\ra |(DK_t)A_0(\varphi_t)|(B(R))\in
L^1([0,T])$. Now we consider the second term on the right hand side
of \eqref{sect-3-prop-1.1}. Again the quantities
$|\det(J_t)|^{-1},\,K_t$ and $J_t^\ast$ are bounded on $[0,T]\times
B(R)$, thus we only need to show that $t\ra
|(\varphi_t^{-1})_\#DA_0|(B(R))$ is integrable on $[0,T]$. We have
  $$|(\varphi_t^{-1})_\#DA_0|(B(R))=|DA_0|(\varphi_t(B(R))).$$
By (1) of Lemma \ref{sect-3-lem-2}, the set $\cup_{0\leq t\leq
T}\varphi_t(B(R))\subset B(F(1+R^\alpha))$ is bounded, from this we
conclude that $[0,T]\ni t\ra |(\varphi_t^{-1})_\#DA_0|(B(R))$ is a
bounded function, hence integrable.

Now we check the last condition in (2) of Theorem
\ref{sect-2-thm-2}. By Lemma \ref{sect-2-lem-1} (2), we have
  $$|\div(\tilde A_0(t))|\leq |\<\div(K_t),A_0(\varphi_t)\>|
  +|\div(A_0)\circ\varphi_t|.$$
Since $\div(K_t)(x)$ and $A_0(\varphi_t(x))$ are bounded on the
product $[0,T]\times B(R)$, we know that $t\ra
\|\<\div(K_t),A_0(\varphi_t)\>\|_{L^\infty(B(R))}$ is a integrable
function on $[0,T]$. By the local boundedness of $\div(A_0)$ and
Lemma \ref{sect-3-lem-2} (1), we obtain the integrability of
$[0,T]\ni t\ra\|\div(A_0)\circ\varphi_t\|_{L^\infty(B(R))}$. Summing
up these discussions, we arrive at
  $$\int_0^T\|\div(\tilde A_0(t))\|_{L^\infty(B(R))}\d t<+\infty.$$
Therefore all the conditions of Theorem \ref{sect-2-thm-2} are
verified, and we complete the proof. \fin

\medskip

Now we are at the position to give the

\medskip

\noindent{\bf Proof of Theorem \ref{sect-1-thm-1}.} (Existence) We
only need to show that the flow $X_t(x)=\varphi_t(Y_t(x))$ satisfies
the Stratonovich SDE \eqref{SDE}. Remark that for all $w\in\Omega_0$
and any $t\in[0,T]$, the map $X_t:\R^d\ra\R^d$ is well defined
almost everywhere. By the generalized It\^{o} formula (see Theorem
3.3.2 in \cite{Kunita90}) and the definitions of $\varphi_t$ and
$Y_t$, we have for a.e. $x\in\R^d$,
  \begin{eqnarray*}
  \d X_t(x)&=&(\d\varphi_t)(Y_t(x))+(\partial_x\varphi_t)(Y_t(x))\,\d
  Y_t(x)\cr
  &=&\bigg(\sum_{i=1}^mA_i(\varphi_t)\circ\d w^i_t\bigg)(Y_t(x))
  +J_t(Y_t(x))K_t(Y_t(x))A_0\big[\varphi_t(Y_t(x))\big]\d t\cr
  &=&\sum_{i=1}^mA_i(X_t(x))\circ\d w^i_t+A_0(X_t(x))\,\d t.
  \end{eqnarray*}
To show the quasi-invariance of the flow $X_t$, let $\tilde\rho_t$
be the Radon-Nikodym density of $(Y_t)_\#\L_d$ with respect to
$\L_d$, then for any $\phi\in C^\infty_c(\R^d)$,
  \begin{eqnarray*}
  \int_{\R^d}\phi(X_t(x))\,\d x&=&\int_{\R^d}\phi\big[\varphi_t(Y_t(x))\big]\,\d
  x\cr
  &=&\int_{\R^d}\phi[\varphi_t(y)]\tilde\rho_t(y)\,\d y
  =\int_{\R^d}\phi(x)\big(\tilde\rho_t|\det(K_t)|\big)(\varphi_t^{-1}(x))\,\d
  x,
  \end{eqnarray*}
where the last equality follows from the change of variable. Hence
  $$(X_t)_\#\L_d=\big(\tilde\rho_t|\det(K_t)|\big)(\varphi_t^{-1})\,\L_d.$$

(Uniqueness) Suppose there is another solution $Z_t$, we consider
$\tilde Z_t=\varphi_t^{-1}(Z_t)$. We will show that $\tilde Z_t$
solves the ODE \eqref{sect-3.2}. In fact, by \cite{Bismut} (see pp.
103--106, or (5.1) in \cite{FangLuo07}),
  \begin{equation}\label{sect-3.4}
  \d\varphi_t^{-1}(x)=-K_t\big(\varphi_t^{-1}(x)\big)
  \bigg(\sum_{i=1}^m A_i(x)\circ\d w^i_t\bigg).
  \end{equation}
Again by the generalized It\^{o} formula (\cite{Kunita90} Theorem
3.3.2),
  $$\d\tilde
  Z_t=(\d\varphi_t^{-1})(Z_t)+\big[(\partial_x\varphi_t^{-1})(Z_t)\big]
  \circ\d Z_t.$$
Recall that for a.e. $x\in\R^d$, $Z_t(x)$ solves the Stratonovich
SDE \eqref{SDE}, and $\partial_x\varphi_t^{-1}=
K_t\big(\varphi_t^{-1}(x)\big)$. Combining these results with
\eqref{sect-3.4} gives rise to
  \begin{eqnarray}
  \d\tilde Z_t(x)&=&-K_t\big(\tilde Z_t(x)\big)
  \bigg(\sum_{i=1}^mA_i(Z_t)\circ\d w^i_t\bigg)\cr
  &&+K_t\big(\varphi_t^{-1}(Z_t(x))\big)\circ\bigg(\sum_{i=1}^m
  A_i(Z_t(x))\circ\d w^i_t+A_0(Z_t(x))\,\d t\bigg)\cr
  &=&K_t(\tilde Z_t(x))A_0\big[\varphi_t(\tilde Z_t(x))\big]\d t
  =\tilde A_0(t,\tilde Z_t(x))\,\d t.
  \end{eqnarray}
That is, $\tilde Z_t(x)$ solves the ODE \eqref{sect-3.2} for a.e.
$x\in\R^d$. But by Proposition \ref{sect-3-prop-1}, this equation
generates a unique flow $Y_t:\R^d\ra\R^d$. Hence $Y_t=\tilde
Z_t=\varphi_t^{-1}(Z_t)$ for almost every $x\in\R^d$, which implies
that any solution $Z_t$ to \eqref{SDE} can be expressed as the
composition $\varphi_t(Y_t)$. We get the uniqueness of \eqref{SDE}.
\fin

\medskip

By Proposition \ref{sect-3-prop-1}, we know that the stochastic flow
$X_t$ in Theorem \ref{sect-1-thm-1} has a inverse flow
$X_t^{-1}=Y_t^{-1}\circ\varphi_t^{-1}$ which consists of measurable
maps on $\R^d$.

\section{Stability of \eqref{SDE} and stochastic transport equation}

In this section we study the stability of the Stratonovich SDE,
proving Theorem \ref{sect-1-thm-2}. As an application we will give
an explicit solution to the corresponding stochastic transport
equation.

\begin{proposition}\label{sect-4-prop-1}
Suppose the conditions of Theorem \ref{sect-1-thm-2}. For each
$n\geq 1$, define $\tilde A^n_0(t,x):=K_t(x)A^n_0(\varphi_t(x))$.
Let $Y^n_t$ be the unique flow generated by the ODE \eqref{sect-3.2}
with $\tilde A_0$ being replaced by $\tilde A^n_0$. Then almost
surely, for any $T,\,R>0$, we have
  $$\lim_{n\ra\infty}\int_{B(R)}\sup_{0\leq t\leq T}|Y^n_t(x)- Y_t(x)|\,\d x=0.$$
\end{proposition}

\Proof. As in Proposition \ref{sect-3-prop-1}, now we check the
conditions in \cite{Ambrosio04} Theorem 6.6. It is clear that
condition (6.3) is verified. Remark again that the uniform
boundedness assumption in (6.4) can be relaxed to allow uniform
linear growth. Since the vector fields $A^n_0$ have the same
sublinear growth, we can prove a uniform growth estimate for $\tilde
A^n_0(t,x)$ similar to the one given in (2) of Lemma
\ref{sect-3-lem-2} (see \cite{FangLuo07} Lemma 2.2 for a proof).
Next we check that $\tilde A^n_0$ converges in
$L^1_{loc}((0,T)\times\R^d, \R^d)$ to $\tilde A_0$ defined in
\eqref{sect-3.1}. We have
  $$|\tilde A^n_0(t,x)-\tilde A_0(t,x)|\leq \|K_t(x)\|\cdot
  |A^n_0(\varphi_t(x))-A_0(\varphi_t(x))|.$$
Since $(t,x)\ra K_t(x)$ is continuous on $[0,T]\times B(R)$, there
is $C_{T,R}>0$ such that
  $$\sup\{\|K_t(x)\|:(t,x)\in[0,T]\times B(R)\}\leq C_{T,R}.$$
Therefore
  \begin{eqnarray*}
  \int_0^T\!\!\!\int_{B(R)}|\tilde A^n_0(t,x)-\tilde A_0(t,x)|\,\d x\d t
  &\leq&C_{T,R}\int_0^T\!\!\!\int_{B(R)}
  |A^n_0(\varphi_t(x))-A_0(\varphi_t(x))|\,\d x\d t\cr
  &=&C_{T,R}\int_0^T\!\!\!\int_{\varphi_t(B(R))}|A^n_0-A_0|
  \cdot|\det(K_t(\varphi_t^{-1}))|\,\d x\d t.
  \end{eqnarray*}
By Lemma \ref{sect-3-lem-2}, the set $\cup_{0\leq t\leq
T}\varphi_t(B(R))\subset B\big(F(1+R^\alpha)\big)$ is bounded, and
the function $|\det(K_t(x))|$ is bounded on $[0,T]\times B(R)$. As a
consequence,
  $$\int_0^T\!\!\!\int_{B(R)}|\tilde A^n_0(t,x)-\tilde A_0(t,x)|\,\d x
  \leq C'_{T,R}T\int_{B(F(1+R^\alpha))}|A^n_0-A_0|\,\d x\d t,$$
which, by condition (2) of Theorem \ref{sect-1-thm-2}, converges to
0 for almost surely $w\in\Omega_0$.

Now we verify the condition (6.5) of Theorem 6.6 in
\cite{Ambrosio04}. By the definition of $\tilde A^n_0(t,x)$, it
holds
  $$\nabla\tilde A^n_0(t,x)=(\nabla K_t(x))A^n_0(\varphi_t(x))
  +K_t(x)(\nabla A^n_0)(\varphi_t(x))J_t(x),$$
hence
  $$\|\nabla\tilde A^n_0(t,x)\|\leq \|\nabla K_t(x)\|\cdot |A^n_0(\varphi_t(x))|
  +\|K_t(x)\|\cdot\|J_t(x)\|\cdot\|(\nabla A^n_0)(\varphi_t(x))\|.$$
The terms $\|\nabla K_t(x)\|,\,\|K_t(x)\|$ and $\|J_t(x)\|$ are
bounded on $[0,T]\times B(R)$. By Lemma \ref{sect-3-lem-2} and the
fact that $A^n_0$ have uniform growth, it is easy to show that the
sequence $|A^n_0(\varphi_t(x))|$ has an upper bound on $[0,T]\times
B(R)$ independent of $n$. Regarding the last term, notice again that
$\cup_{0\leq t\leq T}\varphi_t(B(R))$ is a bounded subset and that
$\|\nabla A^n_0\|$ are locally bounded. Summing up the above
arguments we obtain the boundedness of $\nabla\tilde A^n_0(t,x)$ on
$[0,T]\times B(R)$.

Finally to verify the condition (6.6) in \cite{Ambrosio04}, in view
of Remark 6.3, it is sufficient to show that for all $R>0$
  $$C_T:=\sup_{n\geq 1}\int_0^T\|\div(\tilde
  A^n_0(t,\cdot))\|_{L^\infty(B(R))}\,\d t<+\infty.$$
By the definition of $\tilde A^n_0$ and (2) of Lemma
\ref{sect-2-lem-1}, we have
  $$|\div(\tilde A^n_0(t,\cdot))|\leq |\div(K_t)|\cdot|A^n_0\circ\varphi_t|
    +|\div(A^n_0)\circ\varphi_t|.$$
Similar discussions as above lead to the desired result. Thus all
the conditions in \cite{Ambrosio04} Theorem 6.6 are satisfied, and
we complete the proof. \fin

\begin{corollary}\label{sect-4-cor-1}
Under the conditions of Theorem \ref{sect-1-thm-2}, for any
$p\geq1$, almost surely
  $$\lim_{n\ra\infty}\int_{B(R)}\sup_{0\leq t\leq T}
  |Y^n_t(x)- Y_t(x)|^p\,\d x=0.$$
\end{corollary}

\Proof. By Proposition \ref{sect-4-prop-1}, the sequence
$\sup_{0\leq t\leq T} |Y^n_t(\cdot)- Y_t(\cdot)|$ converges to 0 in
the Lebesgue measure on the ball $B(R)$. Hence we only need to show
that this sequence is bounded in $L^p(B(R),\d x)$ for any $p>1$,
then the desired result follows from the uniform integrability.

By the growth estimate of $\tilde A_0$ in Lemma \ref{sect-3-lem-2},
it is easy to deduce that $|Y_t(x)|\leq \tilde \Phi\cdot(1+|x|)$,
where $\tilde \Phi\in\cap_{p>1}L^p(\Omega)$ (see (iii) in the proof
of \cite{Luo09} Lemma 3.2). Remark that for every $n\geq1$, $\tilde
A^n_0$ has the same growth as $\tilde A_0$, hence
$\sup_{n\geq1}|Y^n_t(x)|\leq \tilde \Phi\cdot(1+|x|)$. Therefore
  \begin{equation}\label{sect-4-cor-1.1}
  \sup_{0\leq t\leq T} |Y^n_t(x)- Y_t(x)|^p
  \leq 2^p\tilde\Phi^p(1+|x|^p),
  \end{equation}
which implies clearly the boundedness of the sequence in
$L^p(B(R),\d x)$. \fin

\medskip

Now we can prove Theorem \ref{sect-1-thm-2}.

\medskip

\noindent{\bf Proof of Theorem \ref{sect-1-thm-2}.} We use the
representations of the solutions: $X_t=\varphi_t(Y_t)$ and $X^n_t=
\varphi_t(Y^n_t)$. By the mean value formula,
  \begin{eqnarray}\label{sect-4.1}
  X^n_t(x)-X_t(x)&=&\varphi_t(Y^n_t(x))-\varphi_t(Y_t(x))\cr
  &=&\bigg(\int_0^1J_t\big((1-u)Y_t(x)+uY^n_t(x)\big)\d
  u\bigg)(Y^n_t(x)-Y_t(x)).
  \end{eqnarray}
By Lemma \ref{sect-3-lem-2} and the growth estimates of $Y_t$ and
$Y^n_t$ given in the proof of Corollary \ref{sect-4-cor-1},
  \begin{eqnarray*}
  \big\|J_t\big((1-u)Y_t(x)+uY^n_t(x)\big)\big\|
  &\leq&G\big(1+|Y_t(x)|^\beta+|Y^n_t(x)|^\beta\big)\cr
  &\leq&G\big(1+2\tilde\Phi^\beta\big)\big(1+|x|^\beta\big).
  \end{eqnarray*}
Therefore by \eqref{sect-4.1},
  \begin{eqnarray}\label{sect-4.2}
  |X^n_t(x)-X_t(x)|&\leq&
  \bigg(\int_0^1\big\|J_t\big((1-u)Y_t(x)+uY^n_t(x)\big)\big\|\d u\bigg)
  |Y^n_t(x)-Y_t(x)|\cr
  &\leq&G\big(1+2\tilde\Phi^\beta\big)\big(1+|x|^\beta\big)|Y^n_t(x)-Y_t(x)|.
  \end{eqnarray}
As a result, for $w\in \Omega_0$,
  \begin{equation*}
  \int_{B(R)}\sup_{0\leq t\leq T}|X^n_t(x)-X_t(x)|\,\d x
  \leq G\big(1+2\tilde\Phi^\beta\big)\big(1+|R|^\beta\big)
  \int_{B(R)}\sup_{0\leq t\leq T}|Y^n_t(x)-Y_t(x)|\,\d x,
  \end{equation*}
which converges to 0 almost surely by Proposition
\ref{sect-4-prop-1}.

Now we prove the $L^p(\Omega)$ convergence of solutions for any
$p\geq1$. Indeed we will prove a stronger result: for any $p\geq1$,
  \begin{equation*}
  \lim_{n\ra\infty}\E\int_{B(R)}\sup_{0\leq t\leq T}
  |X^n_t(x)-X_t(x)|^p\,\d x=0.
  \end{equation*}
Similar to Corollary \ref{sect-4-cor-1}, it is enough to show that
the sequence $\big\{\int_{B(R)} \sup_{0\leq t\leq
T}|X^n_t(x)-X_t(x)|^p\,\d x:n\geq1\big\}$ is bounded in some
$L^q(\Omega)\ (q>1$). However this follows easily from
\eqref{sect-4-cor-1.1}, \eqref{sect-4.2} and the facts that $G,\,
\tilde\Phi$ belong to all $L^q(\Omega)$. The proof of Theorem
\ref{sect-1-thm-2} is complete. \fin

\medskip

As an application of the stability result, now we study the
corresponding stochastic transport equation with the purpose of
constructing an explicit solution to it, by using the flow generated
by \eqref{SDE}. First we give some remarks on the ``inverse'' flow
associated to \eqref{SDE}. To this end we regularize the drift
vector field $A_0$ by convolution using a standard kernel $\chi_n$:
$A^n_0=\phi_n\cdot(A_0\ast \chi_n)$, here $\phi_n(x)=\phi(x/n)$ with
$\phi\in C^\infty_c(\R^d, [0,1])$ satisfying
  $$\phi|_{B(1)}\equiv1,\quad \supp(\phi)\subset B(2).$$
Let $X^n_t(x,w)$ be the smooth flow associated to \eqref{SDE} with
$A_0$ being replaced by $A^n_0$. Then it is clear that the
conditions in Theorem \ref{sect-1-thm-2} are satisfied, hence for
any $p>1$ and $T,R>0$, we have
  \begin{equation}\label{sect-4.3}
  \lim_{n\ra\infty}\E\int_{B(R)}\sup_{0\leq t\leq T}|X^n_t(x)-X_t(x)|^p\,\d x=0.
  \end{equation}
Fix some $T>0$, define the time-reversed Brownian motion $\hat
w^T_t=w_T-w_{T-t}$ and consider
  \begin{equation}\label{inverseSDE}
  \d\hat X^T_t=\sum_{i=1}^mA_i\big(\hat X^T_t\big)\circ\d\big(\hat w^T_t\big)^i
  -A_0\big(\hat X^T_t\big)\,\d t,\quad X_0=x.
  \end{equation}
Under the conditions of Theorem \ref{sect-1-thm-1}, this equation
still generates a unique flow $\hat X^T_t,\,0\leq t\leq T$.
Similarly we have $\hat X^{n,T}_t\big(x,\hat w^T\big)$ which is the
solution of the above equation by replacing $A_0$ with $A^n_0$, then
we have $(X^n_T)^{-1}=\hat X^{n,T}_T$ and
  \begin{equation}\label{sect-4.4}
  \lim_{n\ra\infty}\E\int_{B(R)}\sup_{0\leq t\leq T}
  |\hat X^{n,T}_t(x)-\hat X^T_t(x)|^p\,\d x=0.
  \end{equation}

With the convergence results \eqref{sect-4.3} and \eqref{sect-4.4}
in hand, we can prove

\begin{lemma}\label{sect-4-lem-1}
Assume the conditions of Theorem \ref{sect-1-thm-1} and that the
divergence $\div(A_0)$ is bounded on $\R^d$. Then for every $T>0$,
$X_T^{-1}=\hat X^T_T$ a.e. on $\R^d$, and the density function
$\sigma_T$ of $(X_T^{-1})_\#\L_d$ with respect to $\L_d$ has the
expression:
  $$\sigma_T(x)=\exp\bigg(\sum_{i=1}^m\int_0^T\div(A_i)(X_s(x))\circ\d w^i_s
  +\int_0^T\div(A_0)(X_s(x))\,\d s\bigg).$$
\end{lemma}

\Proof. For every $n\geq1$, we have $((X^n_t)^{-1})_\#\L_d=
\sigma^{(n)}_t\L_d$ where (see Lemma 4.3.1 in \cite{Kunita90})
  $$\sigma^{(n)}_t(x)=\exp\bigg(\sum_{i=1}^m\int_0^t\div(A_i)(X^n_s(x))\circ\d w^i_s
  +\int_0^t\div(A^n_0)(X^n_s(x))\,\d s\bigg).$$
Next $\div(A^n_0)=\<\nabla\phi_n,A_0\ast\chi_n\>+
\phi_n\cdot(\div(A_0)\ast\chi_n)$. Since $A_0$ has sublinear growth,
it is obvious that there exists $C>0$ such that $\sup_{n\geq
1}|(A_0\ast\chi_n)(x)|\leq C(1+|x|)$. By the definition of $\phi_n$,
  \begin{eqnarray*}
  |\<\nabla\phi_n,A_0\ast\chi_n\>|&\leq&
  \frac1n|\nabla\phi(\cdot/n)|\cdot|A_0\ast\chi_n|\cr
  &\leq&\frac {C\|\nabla\phi\|_\infty}n (1+|x|)\ch_{\{n\leq |x|\leq 2n\}}
  \leq 3C\|\nabla\phi\|_\infty.
  \end{eqnarray*}
Hence the divergences
  $$|\div(A^n_0)|\leq |\<\nabla\phi_n,A_0\ast\chi_n\>|
  +|\phi_n\cdot(\div(A_0)\ast\chi_n)|
  \leq 3C\|\nabla\phi\|_\infty+\|\div(A_0)\|_\infty$$
are uniformly bounded on $\R^d$. As a result (see Lemma 3.5 in
\cite{Luo09}), for any $p\in\R$,
  $$\sup_{n\geq 1}\sup_{0\leq t\leq T}\sup_{x\in\R^d}
  \E\big[(\sigma^{(n)}_t(x))^p\big]<+\infty.$$
Now similar arguments as in the proof of Theorem 2.6 in
\cite{Zhang09} lead to the result. \fin

\medskip

This lemma tells us that if $\div(A_0)$ is bounded, then the flow we
constructed in Theorem \ref{sect-1-thm-1} is indeed an almost
everywhere stochastic invertible flow in the sense of Definition 2.1
in \cite{Zhang09}. Furthermore by \eqref{sect-4.4}, for any $t>0$,
it holds
  \begin{equation}\label{sect-4.4.5}
  \lim_{n\ra\infty}\E\int_{B(R)}
  \big|(X^{n}_t)^{-1}(x)-X^{-1}_t(x)\big|^p\,\d x=0.
  \end{equation}

Now we can construct an explicit solution to the corresponding
stochastic transport equation by using the inverse flow $X_t^{-1}$,
as in Theorem 5.1 of \cite{FangLuo07}. Though we follow the idea of
the proof of \cite{FangLuo07} Theorem 5.1, the difference is that
here we only assume that the initial value $\theta_0$ is measurable,
hence the proof of \eqref{sect-4-prop-2.10} is much more delicate
than (5.18) in \cite{FangLuo07} (see \cite{Zhang09} Proposition 2.3
for a different method, but $\theta_0$ is supposed to be bounded
there).

\begin{proposition}\label{sect-4-prop-2}
Assume the conditions in Theorem \ref{sect-1-thm-1} and that
$\div(A_0)$ is bounded. Let $\theta_0:\R^d\ra\R$ be a measurable
function with polynomial growth. Then $\theta(t,x):=\theta_0
(X_t^{-1}(x))$ is a distributional solution to the following
stochastic transport equation:
  \begin{equation}\label{sect-4-prop-2.0}
  \d \theta(t)=-\sum_{i=1}^m\<\nabla \theta(t),A_i\>\circ\d w^i_t
  -\<\nabla \theta(t),A_0\>\,\d t,\quad \theta|_{t=0}=\theta_0,
  \end{equation}
that is, for any $\phi\in C^\infty_c(\R^d)$,
  \begin{equation}\label{sect-4-prop-2.1}
  (\theta(t),\phi)_{L^2}=(\theta_0,\phi)_{L^2}+\sum_{i=1}^m\int_0^t(\theta(s),\div(\phi A_i))_{L^2}
  \circ\d w^i_s+\int_0^t(\theta(s),\div(\phi A_0))_{L^2} \d s.
  \end{equation}
where $(\cdot,\cdot)_{L^2}$ is the inner product in $L^2(\R^d,\d
x)$.
\end{proposition}

\Proof. First we transform the equation \eqref{sect-4-prop-2.1} into
the It\^{o} form:
  \begin{eqnarray}\label{sect-4-prop-2.2}
  (\theta(t),\phi)_{L^2}&=&(\theta_0,\phi)_{L^2}+\sum_{i=1}^m\int_0^t(\theta(s),\div(\phi A_i))_{L^2}
  \d w^i_s+\int_0^t(\theta(s),\div(\phi A_0))_{L^2} \d s\cr
  &&+\frac12\sum_{i=1}^m\int_0^t\big(\theta(s),\div(\div(\phi
  A_i)A_i)\big)_{L^2}\d s.
  \end{eqnarray}
The proof is similar to that of Theorem 5.1 in \cite{FangLuo07} and
we divide it into two parts.

{\bf Step 1.} We assume $\theta_0\in C^\infty_c$. In this step,
similar to the proof of Theorem 5.1 in \cite{FangLuo07}, the key
point is to show
  $$\lim_{n\ra\infty}\E\bigg(\sup_{|x|\leq R}|\theta_n(t,x)
  -\theta(t,x)|^p\bigg)=0,$$
where $\theta_n(t,x)=\theta_0\big((X^n_t)^{-1}(x)\big)$ and $X^n_t$
is defined before Lemma \ref{sect-4-lem-1}. However, similar to
(5.12) on p.1102 of \cite{FangLuo07}, we only need a weaker form of
the convergence result:
  \begin{equation}\label{sect-4-prop-2.3}
  \lim_{n\ra\infty}\E\int_{B(R)}|\theta_n(t,x)
  -\theta(t,x)|^p\,\d x=0.
  \end{equation}
To this end, notice that $|\theta_n(t,x) -\theta(t,x)|\leq
\|\nabla\theta_0\|_\infty|(X^n_t)^{-1}(x) -X_t^{-1}(x)|$, thus by
\eqref{sect-4.4.5}, we still have
  $$\E\int_{B(R)}|\theta_n(t,x)-\theta(t,x)|^p\,\d x
  \leq \|\nabla\theta_0\|_\infty^p\,\E\int_{B(R)}|(X^n_t)^{-1}(x)
  -X_t^{-1}(x)|^p\d x\ra0,$$
as $n\ra\infty$. Hence \eqref{sect-4-prop-2.3} holds and the rest of
the arguments in Step 1 of the proof of Theorem 5.1 in
\cite{FangLuo07} still work here.

{\bf Step 2.} Now suppose that $\theta_0$ is measurable with
polynomial growth. Define $\theta_0^n=\phi_n\, (\theta_0*\chi_n)$.
Then there exist $C>0$ and $q_0\in\Z_+$ which are independent of $n$
such that
  \begin{equation}\label{sect-4-prop-2.4}
  |\theta_0^n(x)|\leq C\, (1+|x|^{q_0}),\quad
  x\in\R^d.
  \end{equation}
Use again the notation $\theta_n(t,x)$ to denote
$\theta_n(t,x)=\theta_0^n(X_t^{-1}(x))$ where $X_t(x)$ is now the
solution to SDE \eqref{SDE}. Then by Step 1, $\theta_n$ satisfies
\eqref{sect-4-prop-2.2}. Now using the SDE \eqref{inverseSDE} and
the moment estimate, we have for any $T>0$ and $t\in[0,T]$,
  $$\E\bigg(\sup_{0\leq s\leq t} |\hat
  X^t_s(x)|^p\bigg)\leq C_{p,T}\, (1+|x|^p).$$
In particular, for each $t\in [0,T]$,
  $$\E\big(|X_t^{-1}(x)|^p\big)\leq C_{p,T}\, (1+|x|^p).$$
By \eqref{sect-4-prop-2.4}, it holds that
  \begin{equation}\label{sect-4-prop-2.4.5}
  \sup_{t\leq T}\E\big(|\theta_n(t,x)|^p\big)+\sup_{t\leq
  T}\E\big(|\theta(t,x)|^p\big)\leq C_{p,T}\,
  (1+|x|^{pq_0}).
  \end{equation}
Therefore for any $p>2$ and $\phi\in C_0^\infty(\R^d)$ with
$\supp(\phi)\subset B(R)$,
  \begin{eqnarray}\label{sect-4-prop-2.5}
  \hskip-8mm\int_0^T \E\big(|(\theta_n(t),\phi)_{L^2}|^p\big)\,\d t
  &\leq& \bigg(\int_{\R^d}|\phi|^q\,
  \d x\bigg)^{p-1}\int_0^T\!\!\!\int_{B(R)}\E(|\theta_n(t,x)|^p)\,\d x\d t\cr
  \hskip-8mm&\leq& C_{p,T}\,\bigg(\int_{\R^d}|\phi|^q\,
  \d x\bigg)^{p-1}\int_{B(R)}(1+|x|^{pq_0})\,\d x<+\infty.
  \end{eqnarray}
where $q$ is the conjugate number of $p$.

Fix some $M>0$. We have
  \begin{eqnarray}\label{sect-4-prop-2.6}
  \E\int_{B(R)}|\theta_n(t,x)-\theta(t,x)|^2\d x&=&
  \E\int_{X_t^{-1}(B(R))}|\theta^n_0(y)-\theta_0(y)|^2\sigma_t(y)\,\d
  y\cr
  &\leq&\E\bigg(\int_{B(M)}+\int_{X_t^{-1}(B(R))\setminus
  B(M)}\bigg)|\theta^n_0(y)-\theta_0(y)|^2\sigma_t(y)\,\d
  y\cr
  &=:&I^n_1(t)+I^n_2(t).
  \end{eqnarray}
Since $\div(A_0)$ is bounded, we have by Lemma 3.5 in \cite{Luo09}
that
  $$\sup_{0\leq t\leq T}\sup_{x\in\R^d}\E(\sigma_t(x))\leq
  C_T<+\infty,$$
hence
  \begin{eqnarray}\label{sect-4-prop-2.7}
  \limsup_{n\ra\infty}I^n_1(t)&=&
  \limsup_{n\ra\infty}\int_{B(M)}|\theta^n_0(y)-\theta_0(y)|^2\E(\sigma_t(y))\,\d
  y\cr
  &\leq& C_T\lim_{n\ra\infty}\int_{B(M)}|\theta^n_0(y)-\theta_0(y)|^2\,\d y=0,
  \end{eqnarray}
where the last equality follows from the convergence of $\theta^n_0$
to $\theta_0$ in $L^2_{loc} (\R^d)$. By the Cauchy inequality,
  \begin{eqnarray}\label{sect-4-prop-2.8}
  I^n_2(t)&\leq&\bigg(\E\int_{X_t^{-1}(B(R))\setminus B(M)}
  |\theta^n_0(y)-\theta_0(y)|^4\sigma_t(y)\,\d y\bigg)^{1/2}
  \bigg(\E\int_{X_t^{-1}(B(R))\setminus B(M)}
  \sigma_t(y)\,\d y\bigg)^{1/2}\cr
  &=:&\big(I^n_{2,1}(t)I^n_{2,2}(t)\big)^{1/2}.
  \end{eqnarray}
We have by \eqref{sect-4-prop-2.4.5},
  \begin{eqnarray}\label{sect-4-prop-2.9}
  I^n_{2,1}(t)&\leq& \E\int_{X_t^{-1}(B(R))}
  |\theta^n_0(y)-\theta_0(y)|^4\sigma_t(y)\,\d y\cr
  &=&\E\int_{B(R)}\big|\theta^n_0(X_t^{-1}(x))-\theta_0(X_t^{-1}(x))\big|^4\,\d
  x\cr
  &\leq& C'\int_{B(R)}(1+|x|^{4q_0})\,\d x<+\infty.
  \end{eqnarray}
Next the function $\sigma_t\ch_{X_t^{-1}(B(R))\setminus B(M)}$ tends
to 0 as $M$ tends to $+\infty$ for $\P\times\L_d$-a.e.
$(w,y)\in\Omega_0\times \R^d$; moreover
$\sigma_t\ch_{X_t^{-1}(B(R))\setminus B(M)}\leq
\sigma_t\ch_{X_t^{-1}(B(R))}$ and
  $$\E\int_{\R^d}\sigma_t\ch_{X_t^{-1}(B(R))}\,\d y
  =\E\int_{\R^d}\ch_{B(R)}\d y=\L_d(B(R))<+\infty.$$
Hence by the dominated convergence theorem, we have
  $$\lim_{M\ra+\infty}I^n_{2,2}(t)=0.$$
This plus \eqref{sect-4-prop-2.8} and \eqref{sect-4-prop-2.9} tells
us that
  $$\lim_{M\ra+\infty} I^n_{2}(t)=0.$$
Therefore by \eqref{sect-4-prop-2.7}, first letting $n$ goes to
$+\infty$ in \eqref{sect-4-prop-2.6}, and then letting $M$ goes to
$\infty$, we obtain
  $$\lim_{n\ra+\infty}\E\int_{B(R)}|\theta_n(t,x)-\theta(t,x)|^2\d x=0.$$
From this we deduce that for any $\phi\in C^\infty_c(\R^d)$ with
$\supp(\phi)\subset B(R)$,
  \begin{equation}\label{sect-4-prop-2.10}
  \E\big[\big((\theta_n(t),\phi)_{L^2}-(\theta(t),\phi)_{L^2}\big)^2\big]
  \leq \bigg(\int_{\R^d}\phi^2\d x\bigg)\E\int_{B(R)}|\theta_n(t,x)-\theta(t,x)|^2\d x
  \ra0
  \end{equation}
as $n\ra\infty$. Now \eqref{sect-4-prop-2.5} and
\eqref{sect-4-prop-2.10} allow us to pass to the limit and the proof
is complete. \fin

\medskip

Now we discuss the connection between the stochastic transport
equation \eqref{sect-4-prop-2.0} and the following transport
equation associated to the random vector field $\tilde A_0$ defined
in \eqref{sect-3.1}:
  \begin{equation}\label{randomTE}
  \d u_t=-\<\nabla u_t,\tilde A_0(t)\>\,\d t,\quad u|_{t=0}=u_0.
  \end{equation}
To this end we first give some preparations. Recall that $\varphi_t$
is the smooth flow defined at the beginning of Section 3. Let
$\tilde\rho_t=\frac{\d ((\varphi_t)_\#\L_d)}{\d \L_d}$ and
$\rho_t=\frac{\d ((\varphi_t^{-1})_\#\L_d)}{\d \L_d}$ be the
Radon-Nikodym densities. We have the following simple equality:
  \begin{equation}\label{randomTE.1}
  \tilde\rho_t(x)=\big[\rho_t\big(\varphi_t^{-1}(x)\big)\big]^{-1}.
  \end{equation}
Indeed, for any $\psi\in C_c^\infty(\R^d)$, we have
  \begin{eqnarray*}
  \int_{\R^d}\psi(x)\,\d x
  &=&\int_{\R^d}\psi\big[\varphi_t\big(\varphi_t^{-1}(x)\big)\big]\,\d x\cr
  &=&\int_{\R^d}\psi[\varphi_t(y)]\rho_t(y)\,\d y
  =\int_{\R^d}\psi(x)\rho_t\big(\varphi_t^{-1}(x)\big)\tilde\rho_t(x)\,\d x,
  \end{eqnarray*}
which leads to \eqref{randomTE.1} due to the arbitrariness of
$\psi\in C_c^\infty(\R^d)$. Moreover by Lemma 4.3.1 in
\cite{Kunita90}, the density $\rho_t$ has an explicit expression:
  \begin{equation}\label{randomTE.2}
  \rho_t(x)=\exp\bigg(\sum_{i=1}^m\int_0^t\div(A_i)(\varphi_s(x))\circ\d w^i_s\bigg).
  \end{equation}

Now we show that the distributional solutions of
\eqref{sect-4-prop-2.0} and \eqref{randomTE} are related to each
other by the smooth flow $\varphi_t$.

\begin{proposition}\label{sect-4-prop-3}
Suppose that $\theta_t$ is a distributional solution to the
stochastic transport equation \eqref{sect-4-prop-2.0}, then almost
surely, $u_t:=\theta_t (\varphi_t)$ solves the transport equation
\eqref{randomTE} with $u|_{t=0}=\theta_0$.
\end{proposition}

\Proof. For any $\psi\in C_c^\infty(\R^d)$, we have
  \begin{equation}\label{sect-4-prop-3.1}
  \int_{\R^d}\theta_t(\varphi_t)\psi\,\d x
  =\int_{\R^d}\theta_t\psi(\varphi_t^{-1})\tilde\rho_t\,\d x.
  \end{equation}
Now we compute the Stratonovich stochastic differentials of
$\psi(\varphi_t^{-1})$ and $\tilde\rho_t$. By \cite{Bismut} (see pp.
103--106, or (5.1) in \cite{FangLuo07}),
  \begin{equation}\label{sect-4-prop-3.1.5}
  \d\varphi_t^{-1}=-K_t(\varphi_t^{-1})\sum_{i=1}^mA_i(x)\circ\d w^i_t.
  \end{equation}
Hence
  \begin{eqnarray*}
  \d\psi(\varphi_t^{-1})=\big\<(\nabla\psi)(\varphi_t^{-1}),\circ\,\d\varphi_t^{-1}\big\>
  =-\sum_{i=1}^m\big\<(\nabla\psi)(\varphi_t^{-1}),K_t(\varphi_t^{-1})A_i\big\>\circ\d
  w^i_t.
  \end{eqnarray*}
Notice that $\nabla(\psi(\varphi_t^{-1}))
=K_t^\ast(\varphi_t^{-1})(\nabla\psi)(\varphi_t^{-1})$, we obtain
  \begin{eqnarray}\label{sect-4-prop-3.2}
  \d\psi(\varphi_t^{-1})=-\sum_{i=1}^m\big\<\nabla(\psi(\varphi_t^{-1})),A_i\big\>\circ\d
  w^i_t.
  \end{eqnarray}
Next we compute $\d\tilde\rho_t$. By \eqref{randomTE.2},
  $$\d\rho_t=\rho_t\sum_{i=1}^m\div(A_i)(\varphi_t)\circ\d w^i_t,$$
hence we deduce from \eqref{sect-4-prop-3.1.5} and the generalized
It\^{o} formula that
  \begin{eqnarray*}
  \d\big[\rho_t(\varphi_t^{-1})\big]&=&(\d\rho_t)(\varphi_t^{-1})
  +\big\<(\nabla\rho_t)(\varphi_t^{-1}),\circ\,\d\varphi_t^{-1}\big\>\cr
  &=&\rho_t(\varphi_t^{-1})\sum_{i=1}^m\div(A_i)\circ\d w^i_t
  -\sum_{i=1}^m\big\<(K_t^\ast\nabla\rho_t)(\varphi_t^{-1}),A_i\big\>\circ\d
  w^i_t.
  \end{eqnarray*}
Using again the It\^{o} formula and by the relation
\eqref{randomTE.1}, we arrive at
  \begin{eqnarray*}
  \d\tilde\rho_t&=&-\big[\rho_t(\varphi_t^{-1})\big]^{-2}\circ\d\big[\rho_t(\varphi_t^{-1})\big]\cr
  &=&-\tilde\rho_t\sum_{i=1}^m\div(A_i)\circ\d w^i_t
  +\tilde\rho_t^2\sum_{i=1}^m\big\<\nabla(\rho_t(\varphi_t^{-1})),A_i\big\>\circ\d
  w^i_t.
  \end{eqnarray*}
Since $\nabla\tilde\rho_t=-\tilde\rho_t^2\nabla
(\rho_t(\varphi_t^{-1}))$, finally we obtain
  \begin{eqnarray}\label{sect-4-prop-3.3}
  \d\tilde\rho_t&=&-\tilde\rho_t\sum_{i=1}^m\div(A_i)\circ\d w^i_t
  -\sum_{i=1}^m\<\nabla\tilde\rho_t,A_i\>\circ\d
  w^i_t
  =-\sum_{i=1}^m\div(\tilde\rho_tA_i)\circ\d w^i_t.
  \end{eqnarray}

Now by the equalities \eqref{sect-4-prop-3.2},
\eqref{sect-4-prop-3.3} and the fact that $\theta_t$ solves the
stochastic transport equation \eqref{sect-4-prop-2.0}, we have
  \begin{eqnarray*}
  \d\big[\theta_t\psi(\varphi_t^{-1})\tilde\rho_t\big]&=&
  \psi(\varphi_t^{-1})\tilde\rho_t\bigg(-\sum_{i=1}^m\<\nabla\theta_t,A_i\>\circ\d w^i_t
  -\<\nabla\theta_t,A_0\>\,\d t\bigg)\cr
  &&-\sum_{i=1}^m\theta_t\tilde\rho_t\<\nabla(\psi(\varphi_t^{-1})),A_i\>\circ\d w^i_t
  -\sum_{i=1}^m\theta_t\psi(\varphi_t^{-1})\div(\tilde\rho_tA_i)\circ\d
  w^i_t.
  \end{eqnarray*}
Since $\div(\psi(\varphi_t^{-1})\tilde\rho_tA_i)=
\psi(\varphi_t^{-1})\div(\tilde\rho_tA_i)+
\tilde\rho_t\<\nabla(\psi(\varphi_t^{-1})), A_i\>$, we arrive at
  \begin{eqnarray*}
  \d\big[\theta_t\psi(\varphi_t^{-1})\tilde\rho_t\big]&=&
  \psi(\varphi_t^{-1})\tilde\rho_t\bigg(-\sum_{i=1}^m\<\nabla\theta_t,A_i\>\circ\d w^i_t
  -\<\nabla\theta_t,A_0\>\,\d t\bigg)\cr
  &&-\sum_{i=1}^m\theta_t\div\big(\psi(\varphi_t^{-1})\tilde\rho_tA_i\big)\circ\d
  w^i_t.
  \end{eqnarray*}
The above equality should be understood in the distributional sense.
More precisely we have obtained
  \begin{eqnarray*}
  \int_{\R^d}\theta_t\psi(\varphi_t^{-1})\tilde\rho_t\,\d
  x&=&\int_{\R^d}\theta_0\psi\,\d x+\sum_{i=1}^m\int_0^t\bigg(\int_{\R^d}
  \theta_s\div\big(\psi(\varphi_s^{-1})\tilde\rho_sA_i\big)\d
  x\bigg)\circ\d w^i_s\cr
  &&+\int_0^t\!\!\!\int_{\R^d}
  \theta_s\div\big(\psi(\varphi_s^{-1})\tilde\rho_sA_0\big)\d
  x\d s\cr
  &&-\sum_{i=1}^m\int_0^t\bigg(\int_{\R^d}
  \theta_s\div(\psi(\varphi_s^{-1})\tilde\rho_sA_i)\d
  x\bigg)\circ\d w^i_s\cr
  &=&\int_{\R^d}\theta_0\psi\,\d x+\int_0^t\!\!\!\int_{\R^d}
  \theta_s\div\big(\psi(\varphi_s^{-1})\tilde\rho_sA_0\big)\d
  x\d s.
  \end{eqnarray*}
As in Proposition \ref{sect-4-prop-2} we denote by
$(\cdot,\cdot)_{L^2}$ the inner product in $L^2(\R^d,\d x)$. By
\eqref{sect-4-prop-3.1} we have
  \begin{equation}\label{sect-4-prop-3.4}
  (\theta_t(\varphi_t),\psi)_{L^2}=(\theta_0,\psi)_{L^2}
  +\int_0^t\!\!\!\int_{\R^d}
  \theta_s\div\big(\psi(\varphi_s^{-1})\tilde\rho_sA_0\big)\d
  x\d s.
  \end{equation}
By the definition of $\tilde\rho_s$,
  \begin{eqnarray*}
  &&\int_{\R^d}\theta_s\div\big(\psi(\varphi_s^{-1})\tilde\rho_sA_0\big)\d
  x\cr
  &&\hskip6mm=\int_{\R^d}\theta_s\tilde\rho_s\big[\big\<(K_s^\ast\nabla\psi)(\varphi_s^{-1}),A_0\big\>
  +\psi(\varphi_s^{-1})\tilde\rho_s^{-1}\<\nabla\tilde\rho_s,A_0\>
  +\psi(\varphi_s^{-1})\div(A_0)\big]\d x\cr
  &&\hskip6mm=\int_{\R^d}\theta_s(\varphi_s)\big[\big\<\nabla\psi,K_s A_0(\varphi_s)\big\>
  +\psi\big\<(\tilde\rho_s^{-1}\nabla\tilde\rho_s)(\varphi_s),A_0(\varphi_s)\big\>
  +\psi\div(A_0)(\varphi_s)\big]\d x.
  \end{eqnarray*}
Lemma \ref{sect-2-lem-2} leads to
  \begin{eqnarray*}
  \div\big(\tilde
  A_0(s)\big)&=&\<\div(K_s),A_0(\varphi_s)\>+\div(A_0)(\varphi_s)\cr
  &=&\big\<(\tilde\rho_s^{-1}\nabla\tilde\rho_s)(\varphi_s),A_0(\varphi_s)\big\>
  +\div(A_0)(\varphi_s),
  \end{eqnarray*}
it follows that
  $$\int_{\R^d}\theta_s\div\big(\psi(\varphi_s^{-1})\tilde\rho_sA_0\big)\d
  x=\int_{\R^d}\theta_s(\varphi_s)\div\big(\psi\tilde A_0(s)\big)\d x.$$
This plus \eqref{sect-4-prop-3.4} gives us
  $$(\theta_t(\varphi_t),\psi)_{L^2}=(\theta_0,\psi)_{L^2}
  +\int_0^t\big(\theta_s(\varphi_s),\div(\psi\tilde A_0(s))\big)_{L^2}\,\d s,$$
which means that almost surely, $\theta_t(\varphi_t)$ is a
distributional solution to the transport equation \eqref{randomTE}
with initial value $\theta_0$. \fin

\begin{remark}
Originally we intended to prove the uniqueness of the solutions to
the stochastic transport equation \eqref{sect-4-prop-2.0} by using
the above proposition. Indeed, if any solution of \eqref{randomTE}
can be represented as $u_t=u_0\big(Y^{-1}_t\big)$, where $Y_t$ is
the flow generated by $\tilde A_0$, then by the above proposition,
we must have $\theta_t(\varphi_t)=\theta_0\big(Y^{-1}_t\big)$ (since
$\theta_t(\varphi_t)|_{t=0}=\theta_0$), which gives us
  $$\theta_t=\theta_0\big[Y^{-1}_t\big(\varphi_t^{-1}\big)\big]
  =\theta_0\big(X_t^{-1}\big).$$
That is to say, any solution of \eqref{sect-4-prop-2.0} is expressed
as the composition of $\theta_0$ and the inverse flow $X_t^{-1}$.
However, since the divergence $\div(\tilde A_0)$ of $\tilde A_0$ is
unbounded, it is difficult to get a meaningful uniqueness result for
the equation \eqref{randomTE}, see \cite{Ambrosio04, Ambrosio08,
DiPernaLions}.
\end{remark}

\section{Approximate differentiability of the flow generated by \eqref{SDE}}

In this section we study the approximate differentiability of the
stochastic flow $X_t$ associated to the Stratonovich SDE \eqref{SDE}
whose drift coefficient $A_0$ belongs to the Sobolev space
$W^{1,1}_{loc}(\R^d,\R^d)$. To this end, we introduce some notations
and results about maximal functions. For any bounded measurable
subset $U\subset\R^d$ with positive Lebesgue measure $\L_d(U)>0$,
define the average of $f\in L_{loc}^1(\R^d)$ on $U$ by
  $$\bint_U f(x)\,\d x=\frac 1{\L_d(U)}\int_U f(x)\,\d x.$$
Then for any $x\in\R^d$ and $R>0$, the local maximal function $M_R
f$ is defined by
  $$M_R f(x)=\sup_{0<r\leq R}\bint_{B(x,r)}|f(y)|\,\d y,$$
where $B(x,r)=\{y\in\R^d:|y-x|\leq r\}$. Here are some results
regarding the maximal function (see \cite{Stein}; for a proof of the
second result, cf. the Appendix of \cite{FangLuoThalmaier}).

\begin{lemma}\label{sect-5-lem-1}
\begin{enumerate}
\item[\rm(1)] For $R,\rho>0$, there are $C_d,\,C_{d,\rho}>0$ such
that for all $f\in L_{loc}^1(\R^d)$, we have
  $$\int_{B(\rho)}M_Rf(x)\,\d x\leq C_{d,\rho}+C_d\int_{B(R+\rho)}|f(x)|\log(2+|f(x)|)\,\d x$$
and for any $\alpha>0$,
  $$\L_d(x\in B(\rho):M_Rf(x)>\alpha)\leq \frac{C_d}\alpha \int_{B(R+\rho)}|f(x)|\,\d x.$$
\item[\rm(2)] Let $f\in W^{1,1}_{loc}(\R^d)$. Then there is $C_d>0$
(independent of $f$) and a negligible set $N\subset\R^d$, such that
for all $x,y\in N^c$ with $|x-y|\leq R$,
  $$|f(x)-f(y)|\leq C_d|x-y|\big((M_R|\nabla f|)(x)+(M_R|\nabla f|)(y)\big).$$
\end{enumerate}
\end{lemma}

We first prove the following result on the approximate
differentiability of the regular Lagrangian flow generated by a
Sobolev vector field $b$. This is an extension of Corollary 2.5 in
\cite{CrippaLellis} to the case where $b$ has linear growth (see
\cite{CrippaLellis} Corollary 3.5 and
\cite{AmbrosioLecumberryManiglia} Remark 3.8 for more general case,
but therein the divergence of $b$ is assumed to be bounded on
$\R^d$).

\begin{proposition}\label{sect-5-prop-1}
Assume that $b\in L^1\big([0,T],W^{1,1}_{loc}(\R^d,\R^d)\big)$
satisfying
\begin{enumerate}
\item[\rm(i)] $\frac{|b_t(x)|}{1+|x|}\in L^\infty([0,T]\times\R^d)$;
\item[\rm(ii)] for any $R>0$, $\int_0^T\|\div(b_t)\|_{L^\infty(B(R))}\d
t<+\infty$;
\item[\rm(iii)] for any $R>0$, $\int_0^T\!\!\int_{B(R)}|\nabla b_t|\log(2+|\nabla b_t|)\,\d x\d
t<+\infty$.
\end{enumerate}
Let $Y_t$ be the regular Lagrangian flow generated by $b$. Then for
any $R>0$ and $\ee>0$, there exists a Borel set $E\subset B(R)$ such
that $\L_d(B(R)\setminus E)<\ee$ and the restriction $Y_t|_E$ is a
Lipschitz map for any $t\in[0,T]$.

In particular, for any $t\in[0,T]$, $Y_t$ is approximately
differentiable $\L_d$-a.e. in $\R^d$.
\end{proposition}

\Proof. We follow the idea of the proof of Corollary 2.5 in
\cite{CrippaLellis}. For $0\leq t\leq T,\, 0<r\leq2R$ and $x\in
B(R)$, define
  $$Q(t,x,r)=\bint_{B(x,r)}\log\bigg(\frac{|Y_t(x)-Y_t(y)|}r+1\bigg)\d y.$$
From Definition \ref{sect-2-def-2}(1), it follows that for a.e. $x$
and for all $r\in(0,2R]$, the map $t\ra Q(t,x,r)$ is Lipschitz and
  \begin{eqnarray}\label{sect-5-prop-1.1}
  \frac{\d Q}{\d t}(t,x,r)&\leq&\bint_{B(x,r)}\bigg|\frac{\d Y_t}{\d t}(x)-\frac{\d Y_t}{\d t}(y)\bigg|
  \cdot\big(|Y_t(x)-Y_t(y)|+r\big)^{-1}\d y\cr
  &=&\bint_{B(x,r)}\frac{|b_t(Y_t(x))-b_t(Y_t(y))|}{|Y_t(x)-Y_t(y)|+r}\,\d
  y.
  \end{eqnarray}
By condition (i) and Gronwall's inequality, it is easy to show that
  $$|Y_t(x)|\leq (1+R)e^{CT},\quad \mbox{for all }x\in B(R),\ 0\leq t\leq T.$$
Therefore for a.e. $x\in B(R)$ and $y\in B(x,r)$, we have
  $$|Y_t(x)-Y_t(y)|\leq |Y_t(x)|+|Y_t(y)|\leq (1+R)e^{CT}+(1+3R)e^{CT}=2(1+2R)e^{CT}=:\tilde R.$$
Since $(Y_t)_\#\L_d\ll\L_d$, we can apply Lemma
\ref{sect-5-lem-1}(2) to get
  $$|b_t(Y_t(x))-b_t(Y_t(y))|\leq C_d|Y_t(x)-Y_t(y)|\cdot\big[(M_{\tilde R}|\nabla b_t|)(Y_t(x))
  +(M_{\tilde R}|\nabla b_t|)(Y_t(y))\big].$$
Substituting this estimate into \eqref{sect-5-prop-1.1} gives us
  \begin{eqnarray*}
  \frac{\d Q}{\d t}(t,x,r)&\leq&\bint_{B(x,r)}C_d\big[(M_{\tilde R}|\nabla b_t|)(Y_t(x))
  +(M_{\tilde R}|\nabla b_t|)(Y_t(y))\big]\d y\cr
  &=&C_d(M_{\tilde R}|\nabla b_t|)(Y_t(x)) +C_d\bint_{B(x,r)}(M_{\tilde R}|\nabla
  b_t|)(Y_t(y))\,\d y.
  \end{eqnarray*}
Integrating with respect to time, we obtain for all $t\in[0,T]$,
  $$Q(t,x,r)\leq \log 2+C_d\int_0^T(M_{\tilde R}|\nabla b_s|)(Y_s(x))\,\d s
  +C_d\int_0^T\!\!\bint_{B(x,r)}(M_{\tilde R}|\nabla
  b_s|)(Y_s(y))\,\d y\d s.$$
Let $\Phi(x)=\int_0^T(M_{\tilde R}|\nabla b_s|)(Y_s(x))\,\d s$, then
by Fubini's theorem,
  $$Q(t,x,r)\leq\log 2+ C_d\Phi(x)+C_d\bint_{B(x,r)}\Phi(y)\,\d y,\quad\mbox{for all }t\in[0,T].$$
Hence by the definition of the maximal function,
  \begin{equation}\label{sect-5-prop-1.2}
  \sup_{0\leq t\leq T}\sup_{0<r\leq2R}Q(t,x,r)\leq \log
  2+C_d\Phi(x)+C_d(M_{2R}\Phi)(x).
  \end{equation}

For $\eta$ sufficiently small, we have
  \begin{eqnarray}\label{sect-5-prop-1.3}
  \hskip-8mm&&\L_d\big(x\in B(R):\log
  2+C_d\Phi(x)+C_d(M_{2R}\Phi)(x)>1/\eta\big)\cr
  \hskip-8mm&&\hskip6mm\leq \L_d\big(x\in B(R):C_d\Phi(x)>1/(3\eta)\big)
  +\L_d\big(x\in B(R):C_d(M_{2R}\Phi)(x)>1/(3\eta)\big).
  \end{eqnarray}
By Chebyshev's inequality,
  $$\L_d\big(x\in B(R):C_d\Phi(x)>1/(3\eta)\big)
  \leq 3\eta C_d\int_{B(R)}\Phi(x)\,\d x.$$
Using Lemma \ref{sect-5-lem-1}(1), we have
  $$\L_d\big(x\in B(R):C_d(M_{2R}\Phi)(x)>1/(3\eta)\big)
  \leq 3\eta C_dC'_d\int_{B(3R)}\Phi(x)\,\d x.$$
Substituting these two estimates into \eqref{sect-5-prop-1.3} and by
the definition of $\Phi(x)$, we obtain
  \begin{eqnarray*}
  I&:=&\L_d\big(x\in B(R):\log
  2+C_d\Phi(x)+C_d(M_{2R}\Phi)(x)>1/\eta\big)\cr
  &\leq& 3\eta C_d(1+C'_d)\int_{B(3R)}\Phi(x)\,\d x\cr
  &=&3\eta C_d(1+C'_d)\int_0^T\!\!\int_{B(3R)}(M_{\tilde R}|\nabla
  b_t|)(Y_t(x))\,\d x\d t.
  \end{eqnarray*}
Using the density $\tilde\rho_t$ of the flow $Y_t$, we get
  $$I\leq 3\eta C_d(1+C'_d)\int_0^T\!\!\int_{Y_t(B(3R))}(M_{\tilde R}|\nabla
  b_t|)(y)\tilde\rho_t(y)\,\d y\d t.$$
In view of the expression of $\tilde\rho_t$ given in Remark
\ref{sect-2-rem-1}, for any $x\in B(3R)$ and $t\in[0,T]$,
  $$\tilde\rho_t(Y_t(x))=\exp\bigg(-\int_0^t\div(b_s)(Y_s(x))\,\d s\bigg)
  \leq \exp\bigg(\int_0^T\|\div(b_s)\|_{L^\infty(B(R_1))}\d s\bigg)=:L,$$
where $R_1=(1+3R)e^{CT}$. Hence by Lemma \ref{sect-5-lem-1}(1),
  \begin{eqnarray*}
  I&\leq&3\eta C_d(1+C'_d)L\int_0^T\!\!\int_{B(R_1)}(M_{\tilde R}|\nabla
  b_t|)(y)\,\d y\d t\cr
  &\leq&3\eta C_d(1+C'_d)L\int_0^T\bigg[C_{d,R_1}+C''_d\int_{B(R_1+\tilde R)}|\nabla b_t|\log(2+|\nabla b_t|)\,\d
  y\bigg]\d t\cr
  &=:&\eta L_1.
  \end{eqnarray*}
Now for any $\ee>0$, let $\eta=\ee/L_1$, then by
\eqref{sect-5-prop-1.2} and the definition of $I$, we have
  \begin{eqnarray*}
  \L_d\bigg(x\in B(R):\sup_{0\leq t\leq
  T}\sup_{0<r\leq2R}Q(t,x,r)>\frac{L_1}\ee\bigg)
  \leq I\leq \frac\ee{L_1}\cdot L_1=\ee.
  \end{eqnarray*}
Let
  $$E=\bigg\{x\in B(R):\sup_{0\leq t\leq
  T}\sup_{0<r\leq2R}Q(t,x,r)\leq \frac{L_1}\ee\bigg\},$$
then $\L_d(B(R)\setminus E)\leq \ee$ and for any $x\in E,\, 0\leq
t\leq T$ and $0<r\leq2R$, the definition of $Q(t,x,r)$ leads to
  \begin{equation}\label{sect-5-prop-1.4}
  \bint_{B(x,r)}\log\bigg(\frac{|Y_t(x)-Y_t(y)|}r+1\bigg)\d y\leq
  \frac{L_1}\ee.
  \end{equation}

Now fix any $x,y\in E$ and let $r=|x-y|$ which is less than $2R$. We
have by the triangular inequality,
  \begin{eqnarray*}
  \log\bigg(\frac{|Y_t(x)-Y_t(y)|}r+1\bigg)&\leq&
  \log\bigg(\frac{|Y_t(x)-Y_t(z)|+|Y_t(z)-Y_t(y)|}r+1\bigg)\cr
  &\leq&\log\bigg(\frac{|Y_t(x)-Y_t(z)|}r+1\bigg)
  +\log\bigg(\frac{|Y_t(z)-Y_t(y)|}r+1\bigg),
  \end{eqnarray*}
therefore by \eqref{sect-5-prop-1.4},
  \begin{eqnarray*}
  &&\log\bigg(\frac{|Y_t(x)-Y_t(y)|}r+1\bigg)\cr
  &&\hskip6mm=\bint_{B(x,r)\cap B(y,r)}
  \log\bigg(\frac{|Y_t(x)-Y_t(y)|}r+1\bigg)\d z\cr
  &&\hskip6mm\leq\bint_{B(x,r)\cap B(y,r)}\bigg[\log\bigg(\frac{|Y_t(x)-Y_t(z)|}r+1\bigg)
  +\log\bigg(\frac{|Y_t(z)-Y_t(y)|}r+1\bigg)\bigg]\d z\cr
  &&\hskip6mm\leq \tilde
  C_d\bint_{B(x,r)}\log\bigg(\frac{|Y_t(x)-Y_t(z)|}r+1\bigg)\d z
  +\tilde
  C_d\bint_{B(y,r)}\log\bigg(\frac{|Y_t(z)-Y_t(y)|}r+1\bigg)\d z\cr
  &&\hskip6mm\leq 2\tilde C_d\cdot \frac{L_1}\ee,
  \end{eqnarray*}
where $\tilde C_d=\L_d(B(x,r))/\L_d(B(x,r)\cap B(y,r))$ only depends
on the dimension $d$. Therefore
  $$|Y_t(x)-Y_t(y)|\leq re^{2\tilde C_dL_1/\ee}=|x-y|e^{2\tilde C_dL_1/\ee}$$
which holds for all $x,y\in E$. Hence $\textup{Lip}(Y_t|_E)\leq
e^{2\tilde C_dL_1/\ee}$. \fin

\medskip

Now we can prove the main result of this section.

\begin{theorem}\label{sect-5-thm-1}
Assume that $A_1,\cdots,A_m\in  C^{3+\delta}_b(\R^d, \R^d)$ and that
$A_0\in W^{1,1}_{loc}(\R^d,\R^d)$ satisfies
\begin{enumerate}
\item[\rm(1)] $A_0$ has sublinear growth;
\item[\rm(2)] $\div(A_0)$ is locally bounded on $\R^d$;
\item[\rm(3)] for any $R>0$, $\int_{B(R)}\|\nabla A_0\|\log(2+\|\nabla A_0\|)\,\d
x<+\infty$.
\end{enumerate}
Then for a.s. $w\in\Omega_0$, for any $R>0$ and $\delta>0$, there
exists a Borel set $E\subset B(R)$ such that $\L_d(B(R)\setminus
E)<\delta$ and the restriction of the flow $X_t$ to $E$ is a
Lipschitz map for any $t\in [0,T]$. In particular, $X_t$ is
approximately differentiable $\L_d$-a.e. in $\R^d$ for any
$t\in[0,T]$.
\end{theorem}

\Proof. Since $X_t=\varphi_t(Y_t)$ and for a.s. $w\in\Omega_0$, the
map $\varphi_t:\R^d\ra\R^d$ is a $C^2$-diffeomorphism on $\R^d$, we
only have to prove the result for the solution $Y_t$ of the random
ODE \eqref{sect-3.2}. Now we check that $\tilde A_0$ satisfies the
conditions given in Proposition \ref{sect-5-prop-1}. First by the
definition of $\tilde A_0(t,\cdot)$ and the conditions on $A_0$, it
is clear that $\tilde A_0(t,\cdot)\in L^1([0,T],
L_{loc}^1(\R^d,\R^d))$. Moreover
  $$\nabla\tilde A_0(t,x)=(\nabla K_t(x))A_0(\varphi_t(x))
  +K_t(x)(\nabla A_0)(\varphi_t(x))J_t(x),$$
hence
  \begin{equation}\label{sect-5-thm-1.1}
  \|\nabla\tilde A_0(t,x)\|\leq \|\nabla K_t(x)\|\cdot |A_0(\varphi_t(x))|
  +\|K_t(x)\|\cdot\|J_t(x)\|\cdot\|(\nabla A_0)(\varphi_t(x))\|.
  \end{equation}
The terms $\|\nabla K_t(x)\|,\,\|K_t(x)\|$ and $\|J_t(x)\|$ are
bounded on $[0,T]\times B(R)$. By Lemma \ref{sect-3-lem-2} and the
fact that $A_0$ has sublinear growth, it is easy to show that
$|A_0(\varphi_t(x))|$ has an upper bound on $[0,T]\times B(R)$. As
for the last term in \eqref{sect-5-thm-1.1}, noticing that
$L:=\cup_{0\leq t\leq T}\varphi_t(B(R))$ is a bounded subset, we
have
  \begin{eqnarray*}
  \int_0^T\!\!\!\int_{B(R)}\|\nabla\tilde A_0(t,x)\|\,\d x\d t
  &\leq&C_{T,R}+C'_{T,R}\int_0^T\!\!\!\int_{B(R)}\|(\nabla A_0)(\varphi_t)\|\,\d x\d
  t\cr
  &=&C_{T,R}+C'_{T,R}\int_0^T\!\!\!\int_{\varphi_t(B(R))}\|\nabla A_0\|\cdot|\det(K_t)(\varphi_t^{-1})|\,\d x\d
  t\cr
  &\leq&C_{T,R}+C''_{T,R}\int_{L}\|\nabla A_0\|\,\d x<+\infty,
  \end{eqnarray*}
where the last inequality follows from the boundedness of
$|\det(K_t)|$ on $[0,T]\times B(R)$. Hence $\tilde A_0\in
L^1([0,T],W^{1,1}_{loc}(\R^d,\R^d))$.

By Lemma \ref{sect-3-lem-2}, the condition (i) in Proposition
\ref{sect-5-prop-1} is easily checked for $\tilde A_0$. The second
condition (ii) can be verified by using the equality in  Lemma
\ref{sect-2-lem-1}(2), as we have done at the end of the proof of
Proposition \ref{sect-3-prop-1}.

Now we check that $\tilde A_0$ satisfies the condition in
Proposition \ref{sect-5-prop-1}(iii). Again by
\eqref{sect-5-thm-1.1} and the above discussions, we have
$\|\nabla\tilde A_0(t)\|\leq C_{T,R}(1+\|(\nabla
A_0)\circ\varphi_t\|)$. Therefore by the simple inequality
$\log(1+s)\leq s$ for all $s\geq0$, we have
  \begin{eqnarray*}
  \log\big(2+\|\nabla\tilde A_0(t)\|\big)&\leq&
  \log\big[(2+C_{T,R})(2+\|(\nabla A_0)\circ\varphi_t\|)\big]\cr
  &\leq&(1+C_{T,R})+\log\big(2+\|(\nabla A_0)\circ\varphi_t\|\big).
  \end{eqnarray*}
As a result,
  \begin{eqnarray}\label{sect-5-thm-1.2}
  \hskip-8mm&&\|\nabla\tilde A_0(t)\|\log\big(2+\|\nabla\tilde
  A_0(t)\|\big)\cr
  \hskip-8mm&&\hskip6mm\leq
  C_{T,R}(1+C_{T,R})\big(1+\|(\nabla A_0)\circ\varphi_t\|\big)\big[1+\log\big(2+\|(\nabla A_0)\circ\varphi_t\|\big)\big]\cr
  \hskip-8mm&&\hskip6mm\leq C_{T,R}(1+C_{T,R})\big[2\big(1+\|(\nabla A_0)\circ\varphi_t\|\big)
  +\|(\nabla A_0)\circ\varphi_t\|\log\big(2+\|(\nabla A_0)\circ\varphi_t\|\big)\big].
  \end{eqnarray}
Again by the fact that $L:=\cup_{0\leq t\leq T}\varphi_t(B(R))$ is
bounded for any $R>0$, we have by the condition (3) that
  \begin{eqnarray*}
  &&\int_0^T\!\!\!\int_{B(R)}\|(\nabla A_0)\circ\varphi_t\|
  \log\big(2+\|(\nabla A_0)\circ\varphi_t\|\big)\d x\d t\cr
  &&\hskip6mm\leq\int_0^T\!\!\!\int_{\varphi_t(B(R))}\|\nabla A_0\|\log\big(2+\|\nabla
  A_0\|\big)|\det(K_t)\circ\varphi_t^{-1}|\,\d y\d t\cr
  &&\hskip6mm \leq C'_{T,R}\,T\int_{L}\|\nabla A_0\|\log\big(2+\|\nabla
  A_0\|\big)\big]\d y<+\infty,
  \end{eqnarray*}
where $C'_{T,R}=\sup\{|\det(K_t(x))|:(t,x)\in[0,T]\times
B(R)\}<+\infty$. This and \eqref{sect-5-thm-1.2} clearly imply that
  $$\int_0^T\!\!\!\int_{B(R)}\|\nabla\tilde A_0(t,x)\|\log\big(2+\|\nabla\tilde
  A_0(t,x)\|\big)\d x\d t<+\infty.$$
The condition (iii) in  Proposition \ref{sect-5-prop-1} is verified
and the proof is complete. \fin

\end{document}